\documentclass{article}

\usepackage{amsmath,amssymb,mathrsfs,amsthm}
\usepackage{hyperref}
\usepackage{array,enumerate}
\usepackage{WataruCommands}
\usepackage[all]{xy}
\usepackage{bm}% used to bald face vectors
\usepackage{comment}

\DeclareMathOperator{\ram}{ram}

\title{A moving lemma for algebraic cycles with modulus and contravariance}

\author{Wataru Kai\thanks{Mathematical Institute, Tohoku University, 6-3 Aramaki Aoba, Sendai 980-8578, Japan. {\tt kaiw@tohoku.ac.jp} \quad \today }
}

\date{}

\begin{document}

\setcounter{tocdepth}{1}
\setcounter{section}{0}

\maketitle
\begin{abstract}
We prove a moving lemma which implies the contravariance of Bloch-Esnault's additive higher Chow group in smooth affine varieties
and Binda-Saito's higher Chow group (taken in the Nisnevich topology) in smooth varieties equipped with effective Cartier divisors.
The new ingredients in the moving method are 
parallel translation {\em with modulus} in the affine space that involves a new integer parameter,
and Noether's normalization lemma over a Dedekind base.
\end{abstract}

\tableofcontents

%%%%%%%%%%%%%%%%%%%%%%%%%%%%%%%
%%%%%%%%%%%%%%%%%%%%%%%%%%%%%%%
\section*{Introduction}

In recent years, the theory of algebraic cycles with modulus has been an attractive subject.
It concerns the behavior of algebraic cycles at boundaries;
more precisely, with how much multiplicity the cycle intersect a chosen effective Cartier divisor
(called the {\em modulus}).
The notion of modulus dates back at least to class field theory.
In algebraic geometry, this concept has been studied since 1952 when Rosenlicht \cite{Rosenlicht1952} introduced the divisor class group relative to a modulus (on a complete nonsingular curve).

\ \ \

The current development was initiated by Bloch and Esnault \cite{BlochEsnault2003}
who introduced the additive higher Chow group 
$\operatorname{TCH}^*(X,n;m)$ (a definitive version by Park \cite{Park2009}),
and it has been a fruitful subject over the last decade.
It is expected to have a relation to the relative $K$-groups
\begin{equation*}
K_*(X\times \bbA ^1,X\times (m+1)\{ 0\} ),
 \end{equation*}
just as Bloch's higher Chow group
$\CH ^i(X,n)$ is related to the $K$-groups:
\begin{equation*}
K_n(X)_{\bbQ }\cong \bigoplus _{i\ge 0}\CH ^i(X,n)_{\bbQ }
 \end{equation*}
for smooth schemes $X$ over a field.
One of the most remarkable results on this object \cite{Rulling2007} is that for the spectrum of a field, the additive higher Chow groups in a ``Milnor range'' $\operatorname{TCH}^n(\Spec (k),n;m)$
are isomorphic to the big de Rham-Witt $\bbW _{ \{ 1,\dots ,m \} }\Omega ^{n-1}_k $.

\ \ \

In the last couple of years, a more general theory of Chow groups with modulus had been anticipated:
The Chow group of zero-cycles with modulus appeared in the work of
Kerz and Saito \cite{KerzSaito2016}
in relation to the class field theory with wild ramification of varieties over finite fields.
Russell \cite{Russell2013} defined a slightly different version earlier and studied its relation to his Albanese variety with modulus.

Binda and Saito \cite{BS2017} then defined the higher Chow group with modulus $\CH ^i(X|D,n)$
for an arbitrary pair of an algebraic scheme $X$ and an effective Cartier divisor $D$.
It is defined as the homology groups of the cycle complex with modulus $z^i(X|D,\bullet )$.
It contains all the groups above (in suitable versions) as particular cases.
It is expected as a cycle-theoretic cohomology theory corresponding to the relative $K$-theory $K_n(X,D)$.

The complex $z^i(X|D,\bullet )$ is believed to produce a sensible cohomology theory at least if $X$ is smooth.
For example, Binda-Saito \cite[\S 4]{BS2017} have shown a zig-zag of quasi-isomorphisms
$z^1(X|D,\bullet )\simeq \ker (\mcal{O}^*(X)\to \mcal{O}^*(D))[1]$ when $X$ is local and factorial.
Comparison maps between the Zariski or Nisnevich hypercohomology of $z^i(X|D,\bullet )$ and $K$-groups and de Rham(-Witt) cohomologies are known; latest results include
\cite{BS2017}
\cite{RullingSaito}
\cite{KrishnaParkCrys} and \cite{Iwasa2017}.

\ \ \

In spite of being a candidate of a nice cohomology theory, it had been unknown if the additive higher Chow group and the higher Chow group with modulus are contravariant for arbitrary morphisms of smooth schemes. In the projective case, this was settled by Krishna and Park \cite{KrishnaParkAdditive, KrishnaParkModule},
but in the general case (e.g.\ affine) the concept of ``modulus'' gets harder to handle.

The aim of this paper is to provide an affirmative answer to this problem at least {\em locally} by proving a moving lemma.
The moving lemma assures the contravariance of the additive higher Chow group in smooth {\em affine} schemes,
and that of the Nisnevich hypercohomology of Binda-Saito's cycle complex
in pairs $(X,D)$ for which $X\setminus D$ is smooth.

Let us explain our results in more detail:

%%%%%%%%%%%%%%%%%%%%%%%%%%%%%%%
%%%%%%%%%%%%%%%%%%%%%%%%%%%%%%%
\subsection*{Moving lemma}

We will often consider pairs $(X,D)$ consisting of an equi-dimensional scheme $X$ over a base field $k$
and an effective Cartier divisor $D$ on it.

For any integer $i\ge 0$, Binda and Saito defined a subcomplex $z^i(X|D,\bullet )$
of Bloch's cycle complex (cubical version);
in particular, elements of $z^i(X|D,p)$ are represented by cycles on $X\times \bbA ^p$
satisfying certain conditions.
When $D=\emptyset $ it reduces to Bloch's higher Chow theory.

The complex is contravariant for flat maps.
The association
\begin{equation*}
(U\xrightarrow{\text{\'{e}tale}}X)
\mapsto
z^i(U|D_U,\bullet )
 \end{equation*}
defines a presheaf on the small \'{e}tale site on $X$.
We will denote this presheaf by $z^i(-|D_-,\bullet )$.

%%%%%%%%%%%%%
%% Theorem %%
\begin{Definition}
For a finite collection $w$ of irreducible locally closed subsets of $X\setminus D$, define a subcomplex
\begin{equation*}
z^i_w(X|D,\bullet )\subset z^i(X|D,\bullet )
 \end{equation*}
by the condition that $V\in z^i(X|D,p)$ is in $z^i_w$ if and only if its support intersects
$W\times F$ properly for every $W\in w$ and face $F$ of $\bbA ^n$. This extends to a subcomplex of presheaves $z^i_w(-|D_-,\bullet )\subset z^i(-|D_-,\bullet )$.

\end{Definition}
%% Theorem ends %%
%%%%%%%%%%%%%%%%%%

%%%%%%%%%%%%%
%% Theorem %%
\begin{Theorem}[Moving Lemma; {\em see} Theorem \ref{Th:full-statement}]
\label{Th:moving-introduction}
Let $X$ be a scheme of finite type over a field,
$D$ an effective Cartier divisor on $X$
and $w $ a finite set of irreducible locally closed subsets of $X\setminus D$.
Assume $X\setminus D$ is smooth.
Then the above inclusion is a quasi-isomorphism in the Nisnevich topology on $X$:
\begin{equation*}
z^i_w(-|D_-,\bullet )\underset{\underset{\mrm{Nis}}{\sim }}{\hookrightarrow }z^i(-|D_-,\bullet ) .
 \end{equation*}

\end{Theorem}
%% Theorem ends %%
%%%%%%%%%%%%%%%%%%

Actually, we prove it in a more general form (Theorem \ref{Th:full-statement})
which we need in a related work \cite{IwasaKai} with Ryomei Iwasa.

%%%%%%%%%%%%%%%%%%%%%%%%%%%%%%%
%%%%%%%%%%%%%%%%%%%%%%%%%%%%%%%
\subsection*{Noether's normalization lemma over a Dedekind base}

Along the proof of Theorem \ref{Th:moving-introduction}, we prove the following version of Noether's normalization lemma which can be of an independent interest.
It is an analog of \cite[Th.10.2.2]{LevineChow}
and the same method has been applied in \cite{SchmidtStrunk2018}
to prove a version of Gabber's presentation theorem over a Dedekind base.

%%%%%%%%%%%%%
%% Theorem %%
\begin{Theorem}[Theorem \ref{Th:Noether}]\label{Th:Noether-introduction}
Let $X\to B$ be an equi-dimensional morphism of relative dimension $n$, with $B$ the spectrum of a Dedekind domain.
Then locally in the Nisnevich topology on $X$ and $B$, there is a finite surjective map
\begin{equation*}
X\to \bbA ^n_B.
 \end{equation*}

\end{Theorem}
%% Theorem ends %%
%%%%%%%%%%%%%%%%%%

This explains the need of the Nisnevich topology
in Theorem \ref{Th:moving-introduction} from the technical side.

%%%%%%%%%%%%%%%%%%%%%%%%%%%%%%%
%%%%%%%%%%%%%%%%%%%%%%%%%%%%%%%
\subsection*{Functoriality of motivic cohomology}

As a direct consequence of Theorem \ref{Th:moving-introduction}, one gets the following functoriality. We explain how it can be deduced in \S \ref{Sec:preliminaries}.

%%%%%%%%%%%%%
%% Theorem %%
\begin{Theorem}\label{Th:functoriality}
Let $(X',D')$, $(X,D)$ be pairs of equi-dimensional schemes and effective Cartier divisors,
and assume $X\setminus D$ is smooth.
Let $f\colon X'\to X$ be a morphism which induces a morphism $D'\to D$ of subschemes.
Then there is a natural map
\begin{equation*}
f^*\colon z^i(-|D_-,\bullet )
\to f_*\left( z^i(-|D'_- ,\bullet ) \right)
 \end{equation*}
in the Nisnevich local derived category of abelian presheaves on $X$.
Consequently there are natural pull-back maps on the hypercohomology groups
\begin{equation*}
f^*\colon \bbH ^n(X_{\Nis }, z^i(-|D_-,\bullet ))
\to 
\bbH ^n(X'_{\Nis }, z^i(-|D'_-,\bullet )).
 \end{equation*}

\end{Theorem}
%% Theorem ends %%
%%%%%%%%%%%%%%%%%%

Unfortunately, it is known that the Nisnevich hypercohomology groups are not the same as the naive homology groups of the cycle complex in general; see e.g.\ \cite[\S 2.1.4]{BS2017}.
We think that the Nisnevich hypercohomology groups are a better invariant, %than the naive homology groups, 
especially because the $K$-theory satisfies the Nisnevich descent. It would be an interesting problem to ask if the Zariski hypercohomology coincides with the Nisnevich.
We also mention that one gets a product structure in the same derived category, where $D$ and $E$ are two divisors on $X$:
\begin{equation*}
z^i(-|D_-,\bullet )\otimes z^j(-|E_-,\bullet )
\to z^{i+j}(-|(D+E)_-,\bullet ) .
 \end{equation*}

A ``projective'' variant of the functoriality, without the need of Nisnevich localization, was proved by Krishna and Park \cite{KrishnaParkModule}.

%%%%%%%%%%%%%%%%%%%%%%%%%%%%%%%
%%%%%%%%%%%%%%%%%%%%%%%%%%%%%%%
\subsection*{Additive higher Chow groups}

Additive higher Chow groups are obtained as a special case of higher Chow groups with modulus:

%%%%%%%%%%%%%
%% Theorem %%
\begin{Definition}
For a scheme $X$ and modulus $m\ge 1$, define
\begin{equation*}
Tz^i(X,\bullet ;m):= z^i(X\times \bbA ^1\mid X\times (m+1)\{ 0\} ,\bullet -1)
 \end{equation*}
and let us call it the additive cycle complex
(there are several notation/indexing conventions).
Its homology groups
\begin{equation*}
\operatorname{TCH}^i(X,n;m)
:= H_n(Tz^i(X,\bullet ;m))
 \end{equation*}
are called {\em additive higher Chow groups} of $X$.

\end{Definition}
%% Theorem ends %%
%%%%%%%%%%%%%%%%%%

The adjective ``additive'' was given because its values are often something additive, such as extensions of $k$-vector spaces as in the result of R{\"u}lling mentioned above.

One can prove analogs of the above results for this situation, and in particular, gets the functoriality of the additive higher Chow group in smooth {\em affine} schemes (without Nisnevich topology): See Theorem \ref{Th:additive-Chow}.

In the main body of the text, we will stick to the case of higher Chow groups with modulus, and make necessary comments for the additive case as needed.

\ \ \

The projective case had been settled by Krishna and Park \cite{KrishnaParkAdditive}.

\subsubsection*{Organization of the paper}

In \S \ref{Sec:preliminaries}, we recall the definition of the cycle complex and define certain subcomplexes of it.
Along the way, we explain how the functoriality (Theorem \ref{Th:functoriality}) follows from the Moving Lemma \ref{Th:moving-introduction}.

\ \ \

Sections \S\S \ref{Sec:Affine-case}-\ref{Sec:General-case} are devoted to the proof of Moving Lemma \ref{Th:moving-introduction}.
Of course, the whole proof follows the strategy exploited for Bloch's higher Chow theory \cite[Part I, Chap.\ II, \S 3.5]{Levine1998Mixed},
which originates from that devised by Chow \cite{Chow1956}.

In \S \ref{Sec:Affine-case}, we treat a special case of Theorem \ref{Th:moving-introduction}, to which the general case will be reduced.
Here we need a slightly new moving construction which could be called {\em parallel translation with modulus}.

In \S \ref{Sec:linear}, we recall the common method of linear projection used in the reduction process.
There is nothing new here.
I have tried to clarify some details with which I had difficulties when learning this technique through the literature.

In \S \ref{Sec:Noether}, we prove Noether's normalization lemma over a Dedekind base (Theorem \ref{Th:Noether-introduction}), another ingredient used in the reduction process.

Finally in \S \ref{Sec:General-case}, we complete the proof of the moving lemma. 

\ \ \

In a supplementary section \S \ref{Sec:simplicial},
we say a few words about the moving lemma for the simplicial variant of the cycle complex.

\paragraph{Acknowledgements}

This paper is based on the author's Ph.D.\ thesis.
I thank my advisors Prof.\ Shuji Saito and Prof.\ Tomohide Terasoma
for encouragements and guidance.
Discussions with Hiroyasu Miyazaki helped me very much. Especially he introduced the literature \cite{LevineChow} to me when I turned very pessimistic.
I am grateful to Sawako for everyday support.
During the work, I was supported by the Program for Leading Graduate Schools, MEXT, Japan, and by Japan Society for the Promotion of Science as a research fellow (JSPS KAKENHI Grant Number 15J02264).
Lastly I thank my parents for, among thousands of other things, making it possible for me to pursue my study in mathematics.

%%%%%%%%%%%%%%%%%%%%%%%%%%%%%%%
%%%%%%%%%%%%%%%%%%%%%%%%%%%%%%%
\section{Preliminaries}\label{Sec:preliminaries}

In this section, we recall Binda-Saito's cycle complex, and the notion of constructible subsets and functions on noetherian topological spaces.
After that we define certain subcomplexes of the cycle complex determined by data of constructible subsets of the variety.

%%%%%%%%%%%%%%%%%%%%%%%%%%%%%%%
%%%%%%%%%%%%%%%%%%%%%%%%%%%%%%%
\subsection{Cycle complex with modulus}

In this article, an {\it algebraic scheme} is a scheme of finite type over a base field.
Algebraic cycles will always be supported on equi-dimensional algebraic schemes.
This includes the case where a scheme is given over a discrete valuation ring and cycles are on its generic fiber.

By convention, the dimension of a topological space is a non-negative integer or $\infty $ if the space is non-empty, and we set $\dim (\emptyset ):= -\infty $ so that the inequality $\dim (X)\le n$ means $X=\emptyset $ whenever $n\le -1$.

%%%%%%%%%%%%%%%%%%%%%%%%%%%%%%%
%%%%%%%%%%%%%%%%%%%%%%%%%%%%%%%
\subsubsection{}

We write $\square ^p:= \Spec (\bbZ [t_1,\dots ,t_p])$ (the $p$-dimensional cube).
Consider the open immersion $\square ^p\subset (\bbP ^1)^p $ and denote by $F_p$ the effective Cartier divisor 
\begin{equation*}
F_p:= \sum _{j=1}^p \{ t_j=\infty \} \subset (\bbP ^1)^p.
 \end{equation*}
When we work over a base ring, we often tacitly take their base change.

Let $X$ be a scheme equipped with an effective Cartier divisor $D$.
A closed subset $V\subset (X\setminus D)\times \square ^p$ is said to satisfy the {\it modulus condition} if the following is true:
Let $\bar{V}^N$ be the normalization of the (reduced) closure $\bar{V}$ of $V$ in $X\times (\bbP ^1)^p$. Then the inequality of Cartier divisors on $\bar{V}^N$ each of which is obtained as the pull-back:
\begin{equation*}
D|_{\bar{V}^N} \le F_p |_{\bar{V}^N}
 \end{equation*}
holds.

%%%%%%%%%%%%%%%%%%%%%%%%%%%%%%%
%%%%%%%%%%%%%%%%%%%%%%%%%%%%%%%
\subsubsection{}

The {\it faces} of $\square ^p $ are the closed subschemes defined by some of the equations $\{ t_j=\varepsilon \} $
with $1\le j\le p$ and $\varepsilon =0,1$.
A closed subset $V\subset (X\setminus D)\times \square ^p$
is said to satisfy the {\it face condition} if for every face $F$ of $\square ^p$, the intersection $V\times _{\square ^p} F$ has codimension in $V$ at least $\codim _{\square ^p}(F)$.

%%%%%%%%%%%%%
%% Theorem %%
\begin{definition}[Cycle complex with modulus \cite{BS2017}]
For a pair $(X,D)$ of an equi-dimensional algebraic scheme and an effective Cartier divisor on it,
let $\underline{z}^i(X|D,p)$ be the group of codimension $i$ algebraic cycles on $(X\setminus D)\times \square ^p$ whose supports satisfy the modulus and the face conditions.
\end{definition}
%% Theorem ends %%
%%%%%%%%%%%%%%%%%%

The groups $\underline{z}^i(X|D,p)$ ($p\in \bbN $) are known to form a cubical abelian group $\underline{z}^i(X|D,\bullet )$ (an observation of Krishna-Park \cite{KrishnaParkAdditive}).
In particular, we can consider the associated (naive) complex, also denoted by the same symbol, by giving the differentials $\sum _{j=1}^p (-1)^j(\partial _{j,0}^*- \partial _{j,1}^*)$.
Here, $\partial _{j,\vare }\colon \square ^{p-1}\to \square ^p $ is the immersion of the face $\{ t_j=\vare \} $: $(t_1,\dots ,t_{p-1})\mapsto (t_1,\dots ,\overset{\ j}{\check{\vare }},t_j,\dots ,t_{p-1})$.
Its associated non-degenerate complex $z^i(X|D,\bullet )$ is defined (in degree $p$) as the quotient of $\underline{z}^i(X|D,p )$ by the images of $\underline{z}^i(X|D,q)$ ($q< p$) by the degeneracy maps $\square ^p\twoheadrightarrow \square ^q$ (called the degenerate part).
An alternative definition of $z^i(X|D,p)$ is as the subgroup $\bigcap _{j=1}^p \ker ( \partial _{j,0}^* ) $. 
The two definitions are isomorphic by the canonical map.
It also gives a splitting of $\underline{z}^i$ into $z^i$ and the degenerate part.

%%%%%%%%%%%%%%%%%%%%%%%%%%%%%%%
%%%%%%%%%%%%%%%%%%%%%%%%%%%%%%%
\subsection{Constructible subsets}

Let $X$ be a noetherian topological space.
The set $\mfr{C}(X)$ of constructible subsets of $X$ is the smallest subset of the power set $\mfr{P}(X)$ of $X$ satisfying
\begin{itemize}
\item Open sets of $X$ are in $\mfr{C}(X)$;
\item $\mfr{C}(X)$ is closed under finite union and intersection, and complement.
 \end{itemize}
The pull-back of a constructible subset by a continuous map of noetherian spaces is again constructible.
It is not hard to see that the constructible subsets are precisely the unions of finitely many locally closed subsets.

A dense constructible subset of an irreducible noetherian space contains a dense open subset \cite[$0_{\mrm{III}}$ 9.2.2]{EGA}.
Note in particular that a dense constructible subset $Z$ of an algebraic scheme $X$ has the same dimension as the scheme, because $Z$ contains a dense open subset $U$ of $X$ and
therefore $\dim (U)\le \dim (Z)\le \dim (X)$.

%%%%%%%%%%%%%%%%%%%%%%%%%%%%%%%
%%%%%%%%%%%%%%%%%%%%%%%%%%%%%%%
\subsubsection{}

When $T$ is a set, a function $h \colon X\to T$ is said to be {\it constructible}
if its image is a finite set, and for every $t\in T$ the inverse image $h \inv (t)$ is a constructible subset of $X$.
Note in this case that for any subset $S$ of $T$,
its inverse image $h \inv (S)$ is a constructible subset of $X$.

Constructibility of a function $h $ is known to be equivalent to the following condition (\cite[$0_{\mrm{III}}$ 9.3.2]{EGA}; we state it when $X$ is the underlying space of a scheme): For every $x\in X $, there is a non-empty open set of $\overline{\{ x\} }$
where $h $ takes the constant value $h (x)$.

If two functions $h_1,h_2\colon X\rightrightarrows \bbN $ are constructible, then their maximum
$x\mapsto \max \{ h_1(x),h_2(x)\} $ is again constructible.

%%%%%%%%%%%%%
%% Theorem %%
\begin{theorem}[{\cite[IV${}_1$ 1.8.4]{EGA}}]
For every morphism of finite type between noetherian schemes,
the image of a constructible subset is again constructible. 
\end{theorem}
%% Theorem ends %%
%%%%%%%%%%%%%%%%%%

%%%%%%%%%%%%%
%% Theorem %%
\begin{theorem}[{\cite[IV${}_3$ 13.1.3]{EGA}}]\label{Th:Chevalley-dim}
Let $f\colon X\to Y $ be a morphisms of finite type between noetherian schemes.
Then the function
\begin{align*}
X&\to \bbN \\{}
x&\mapsto \dim _x(f\inv (f(x)))
 \end{align*}
is upper semi-continuous; namely, for every interger $n$, the subset $\{ x\in X\mid \dim _s(f\inv (f(x))) \ge n \} $ is closed in $X$.

\end{theorem}
%% Theorem ends %%
%%%%%%%%%%%%%%%%%%

Here, 
the local dimension of a topological space $X$ at a point $x$ is the infimum of the dimensions of its neighborhoods $U$ in $X$:
\begin{equation*}
\dim _x(X) := \inf _{x\in U\subset X} \dim (U).
 \end{equation*}
For example, if $X$ is the union $S\cup C$ of two irreducible components which have dimensions $2$ and $1$ respectively, then the local dimension is $2$ at the points of $S$ whereas it is $1$ at the points of $C\setminus S$.

%%%%%%%%%%%%%
%% Theorem %%
\begin{corollary}\label{Cor:Chevalley-dim}

Let $f\colon X\to Y$ be a morphism of finite type between noetherian schemes.
Let $Z\subset X$ be a constructible subset.
Then the function which 
to $y\in Y$ assigns 
$\dim (f\inv (y)\cap Z)$
is a constructible function on $Y$.

\end{corollary}
%% Theorem ends %%
%%%%%%%%%%%%%%%%%%

%%% Proof begins %%%
\begin{proof}
If $Z=X$, this is a direct consequence of the previous two theorems.

In general, write $Z$ as the union of finitely many locally closed subsets
$Z=\bigcup _i Z_i $, each of which has a scheme structure.
Then we know
$\dim (f\inv (y)\cap Z)=\max _{i} \{ \dim (f\inv (y)\cap Z_i) \} $
because in an algebraic scheme (like $f\inv (y)$),
the dimension of the union of finitely many constructible subsets equals the maximum of their dimensions.
(The dimension is equal to the maximum of the transcendence degrees of the points involved, after all.)

Since each member of the right hand side is already shown to be constructible, our function is constructible.
\end{proof}
%%% Proof ends %%%

%%%%%%%%%%%%%%%%%%%%%%%%%%%%%%%
%%%%%%%%%%%%%%%%%%%%%%%%%%%%%%%
\subsection{Subcomplex determined by constructible subsets}

The following is an essentially equivalent variant of the subcomplex denoted by the symbol $z^i_{\mcal{W},e}$ e.g.\ in \cite{KrishnaParkAdditive, Levine1998Mixed}.
I'd like to try the following new convention partly because this makes it easier to write down the functoriality conditions, and partly because the use of constructible subsets eliminates the choices of decompositions of constructible subsets into locally closed subsets.

%%%%%%%%%%%%%
%% Theorem %%
\begin{definition}\label{Def:increasing-family}
(1)
By an {\it increasing family} $\mcal{C}=\{ C_d \} _{d\in \bbZ }$ of constructible subsets of an equi-dimensional algebraic scheme $X$ we mean a sequence
\begin{equation*}
\cdots \subset C_{0}\subset C_1\subset C_2\subset \cdots 
 \end{equation*}
of constructible subsets where $\dim (C_d)\le d$ and $C_d =X $ for $d\ge \dim (X)$.
In particular, when $d$ is negative, then $C_d$ is the empty set.
We shall drop the adjective {\it increasing} because we are only interested in this type of families.
The family $\mcal{C}_{\mrm{triv}}$ on $X$
defined by
$\emptyset = C_0=\dots =C_{\dim (X)-1} \subsetneq C_{\dim (X)}=X $ will be called the {\it trivial family}.

(2)
When $D\subset X$ is an effective Cartier divisor and $\mcal{C}=\{ C_d\} _{d\in \bbZ } $ is a family of constructible subsets of $X\setminus D$, we define the cubical subgroup 
\begin{equation*}
\underline{z}^i_{\mcal{C}}(\bullet ):= \underline{z}^i_{\mcal{C}}(X|D,\bullet )
 \end{equation*}
of $\underline{z}^i(\bullet ):= \underline{z}^i(X|D,\bullet )$ by the condition: A cycle 
$V\in \underline{z}^i(p)$ is in $\underline{z}^i_{\mcal{C}}(p)$ if and only if for every face $F$ of $\square ^p $ and $d\in \bbZ $, the following inequality of dimensions holds:
\begin{equation*}
\dim ((C_d\times F)\cap |V|)\le d+\dim (F)-i.
 \end{equation*}
Note that for the trivial family $\mcal{C}_{\mrm{triv}}$ on $X\setminus D$, one has
$\underline{z}^i_{\mcal{C}_{\mrm{triv}}}(\bullet )=\underline{z}^i(\bullet )$.

\end{definition}
%% Theorem ends %%
%%%%%%%%%%%%%%%%%%

Given a finite set $w$ of irreducible locally closed subsets of $X\setminus D$ as in Introduction, we define the corresponding $\mcal{C}$ by
\begin{equation*}
C_d:= \bigcup _{\begin{subarray}{c}W\in w \\ \dim (W)\le d\end{subarray}} W   \end{equation*}
so that one recovers the complex $z^i_w(X|D,\bullet )$ in Introduction as $z^i_\mcal{C}(X|D,\bullet )$.

\ \ \

Thus Theorem \ref{Th:moving-introduction}
can be stated in a slightly stronger from that the inclusion
\begin{equation}\label{Eq:principal-case}
\underline{z}^{i}_{\mcal{C}}(-|D_-,\bullet )
\to \underline{z}^{i}(-|D_-,\bullet )
 \end{equation}
be a quasi-isomorphism of presheaves on $X_\Nis $.
Here we regard the left hand side by pulling back constructible subsets by the given {\'e}tale maps $U\to X$.

See Theorem \ref{Th:full-statement} below for the full statement that we prove in this paper.

\begin{comment}
The following is (a special but the principal case of) the main result of the paper, and its proof will be completed at the end of \S \ref{Sec:General-case}.
See Theorem \ref{Th:full-statement} below for the full statement.

%%%%%%%%%%%%%
%% Theorem %%
\begin{theorem}\label{Th:principal-case}
Let $(X,D)$ be a pair of an equi-dimensional $k$-scheme over a field $k$ and an effective Cartier divisor, and $\mcal{C}$ an increasing family {\em (Definition \ref{Def:increasing-family})} of constructible subsets of $X\setminus D$.
Assume $X\setminus D$ is smooth.
%
Then the inclusion map:
%%
\begin{equation*}
\underline{z}^{i}_{\mcal{C}}(X|D,\bullet )
\to \underline{z}^{i}(X|D,\bullet )
 \end{equation*}
%%
is a quasi-isomorphism in the Nisnevich topology of $X$.

\end{theorem}
%% Theorem ends %%
%%%%%%%%%%%%%%%%%%

Here we regard both hand sides as presheaves on the small Nisnevich site of $X$ in an obvious way which is explicated in \S \ref{Sec:def-of-sheaves} below.
%Given an \'{e}tale scheme $U\xrightarrow{\phi }X$, give $U$ the divisor $\phi ^*D$ and the family $\phi ^*\mcal{C}$ of constructible subsets
%$(\phi ^*\mcal{C})_d := \phi \inv (C_d)$.
%Then we can consider the presheaf
%$(U\xrightarrow{\phi }X)\mapsto \underline{z}^i_{\phi ^*\mcal{C}}(U|\phi ^*D,\bullet )$.

\end{comment}

%%%%%%%%%%%%%%%%%%%%%%%%%%%%%%%
%%%%%%%%%%%%%%%%%%%%%%%%%%%%%%%
\subsubsection{}

To illustrate why the quasi-isomorphism \eqref{Eq:principal-case} is useful,
we would like to discuss some functorial behavior of families of constructible subsets.
Let $f\colon X'\to X$ be a map of algebraic schemes and $\mcal{C}'$ be a family of constructible subsets of $X'$.
We define its push-forward $f_*\mcal{C}'$ by
\begin{equation*}
(f_*\mcal{C}')_d :=\bigcup _{e\in \bbZ }
\{ x\in X \mid \dim (C_{d+e}'\times _X x) \ge e
\} \qquad \text{ for } d< \dim (X) 
\end{equation*}
and $(f_*\mcal{C}')_d:= X$ for $d\ge \dim (X)$.
Of course, the union over $0\le e\le \dim (X')$ gives the same subset, and hence it is a constructible subset of $X$ by Proposition \ref{Cor:Chevalley-dim} (Chevalley's theorem).
We obviously have
$(f_*\mcal{C}')_d\subset (f_*\mcal{C}')_{d+1}$ (because $C_{d+e}'\subset C_{d+1+e}'$)
and $\dim ((f_*\mcal{C}')_{d})\le d$ because if we write ${C}_d^{ (e)}$ for the $e$-th summand in the right hand side of the formula, then one should have
\begin{align*}
\dim ({C}_d^{ (e)})+e
&\le
\dim (C_{ d+e}') \\
&\le d+e.
 \end{align*}
When $X$ is also given a family $\mcal{C}=\{ C_d\} _{d\in \bbZ }$ of constructible subsets, we say the morphism $f$ is {\it compatible} with $\mcal{C}'$ and $\mcal{C}$ if the relation
$(f_*\mcal{C}')_d \subset C_d$ holds for all $d\in \bbZ $.

%%%%%%%%%%%%%%%%%%%%%%%%%%%%%%%
%%%%%%%%%%%%%%%%%%%%%%%%%%%%%%%
\subsubsection{ }

Let $(X',D')$ and $(X,D)$ be pairs of algebraic schemes and effective Cartier divisors together with families $\mcal{C}'$ and $\mcal{C}$ of constructible subsets
on $X'\setminus D'$ and $X\setminus D$ respectively. %Assume $Y\setminus E$ is regular so that we can do cycle theory.

Given a {\it morphism of pairs} $f\colon (X',D')\to (X,D)$, i.e., a morphism $X'\to X$ which induces a morphism on the subschemes $D'\to D$,
we may consider the push-forward $f_*\mcal{C}'$ on $X\setminus D$ by seeing the members of $\mcal{C}'$ as constructible subsets of $X'\setminus f\inv (D)$. Consequently, we may ask if $f$ is compatible with the families $\mcal{C}'$ and $\mcal{C}$.
Suppose it is. 
Then if moreover $f$ is of finite Tor dimension (which is a standard assumption when one wants to pull-back algebraic cycles), we have a well-defined pull-back map
$f^*\colon \underline{z}^i_{\mcal{C}}(X|D,\bullet )
\to \underline{z}^i_{\mcal{C}'}(X|D,\bullet )$ as follows.

\ \ \

By the compatibility condition $(f_*\mcal{C}')_{d}\subset C_{d}$, for every $V\in \underline{z}^i_{\mcal{C}}(X|D,p)$ we have:
\begin{align*}
\dim ((C_d' \times F)\cap f\inv (|V|))
&\le \max _{e\in \bbZ }\{ e+\dim ((C_{d-e}\times F)\cap |V|) \}
\\
&\le \max _{e\in \bbZ }\{e+(d-e+\dim (F))-i \}
\\
&= d+\dim (F)-i
 \end{align*}
so that we can define a cycle $f^*V$ on $X'\times \square ^p$ by Serre's Tor formula \cite[Chap.V.C)]{Serre} and it satisfies the necessary dimension conditions for being an element of $\underline{z}^i_{\mcal{C}'}(X'|D',p)$.
Its support satisfies the modulus condition as well because $|V| $ does (cf.\ similarly to the ``containment lemma'' \cite[Prop.2.4]{KrishnaParkAdditive}).
Therefore the pull-back $f^*V$ is in $\underline{z}^i_{\mcal{C}'}(X'|D',p)$.

In particular, by setting $\mcal{C}'$ to be the trivial family and $\mcal{C}:=f_*\mcal{C}'_{\mrm{triv}}$, we get a diagram of inclusion and pull-back 
\begin{equation*}
\underline{z}^i(X|D,\bullet )\supset
\underline{z}^i_{\mcal{C}}(X|D,\bullet )
\xrightarrow{f^*}
\underline{z}^i(X'|D',\bullet ) .
 \end{equation*}
The quasi-isomorphism \eqref{Eq:principal-case} allows us to regard this diagram as a defining datum of a morphism $\underline{z}^i(-|D_-,\bullet )\to f_*(\underline{z}^i(-|D'_-,\bullet ))$ in the derived category $D(X_{\Nis })$.
%See \S \ref{Sec:def-of-sheaves} for the notation.

%%%%%%%%%%%%%%%%%%%%%%%%%%%%%%%
%%%%%%%%%%%%%%%%%%%%%%%%%%%%%%%
\subsubsection{}

As an aside,
the association $((X,D),\mcal{C})\mapsto \underline{z}^i_{\mcal{C}}(X|D,\bullet )$
is naturally a presheaf on the category of pairs $(X,D)$ whose open part $X\setminus D$ is regular and equipped with a family $\mcal{C}$ of constructible subsets with the notion of compatible morphisms above.
We have the forgetful functor $((X,D),\mcal{C})\mapsto (X,D)$ to the category of such pairs, but usually there is no functor in the opposite direction unless we restrict to a subcategory of pairs where every morphism is equi-dimensional
or to special types of diagrams (e.g.\ finite) in the category of pairs.
Given the quasi-isomorphism \eqref{Eq:principal-case},
it is possible to construct a presheaf on the category of pairs $(X,D)$ with $X\setminus D$ smooth which is objectwise quasi-isomorphic to our cycle complex Nisnevich locally, by a categorical construction \cite{Levine2008Coniveau} involving huge homotopy limit and colimit.

%%%%%%%%%%%%%%%%%%%%%%%%%%%%%%%
%%%%%%%%%%%%%%%%%%%%%%%%%%%%%%%
\subsubsection{}\label{Sec:def-of-sheaves}

Next, let $f\colon X'\to X$ be an equi-dimensional morphism of relative dimension $\dim (X'/X)$ inducing a morphism of effective Cartier divisors
$D'\to D$, and suppose that $X$ is given a family $\mcal{C}$ of constructible subsets.
Then we can define its pull-back $f^*\mcal{C}$ by
\begin{equation*}
(f^*\mcal{C})_d:= f\inv (C_{d-\dim (X/Y)}) \setminus f\inv (D).
 \end{equation*}
If moreover $f$ is of finite Tor dimension over $X\setminus D$ (e.g.\ if $f$ is flat or $X\setminus D$ is regular), then we have a natural pull-back map
\begin{equation*}
\underline{z}^i_{\mcal{C}}(X|D,\bullet )
\to \underline{z}^i_{f^*\mcal{C}}(X'|D',\bullet ) .
 \end{equation*}
In particular, for a pair $(X,D)$ equipped with a family $\mcal{C}$, 
we can consider the presheaf 
\begin{equation*}
\underline{z}^{i}_{\mcal{C}}(-|D_-,\bullet )\colon
\quad
(U\xrightarrow{\phi }X)\mapsto \underline{z}^i_{\phi  ^*\mcal{C}}(U|\phi ^*D,\bullet )\end{equation*}
on $X_\Nis $.
When $\mcal{C}$ is the trivial family $\mcal{C}_{\mrm{triv}}$, we will drop it from the notation.

\ \ \

We extend the definition slightly. Given an equi-dimensional scheme $\psi \colon Y\to \Spec (k)$,
we consider the presheaf 
\begin{equation*}
\underline{z}^{i}_{\mcal{C}}(-\times Y\mid D_-\times Y,\bullet )\colon 
\quad
(U\xrightarrow{\phi }X)\mapsto \underline{z}^i_{(\phi \times \psi ) ^*\mcal{C}}(U\times Y\mid \phi ^*D\times Y,\bullet ) 
\end{equation*}
on $X_\Nis $.
Let us denote the value also by
$\underline{z}^i_{\mcal{C}}(U\times Y\mid \phi ^*D\times Y,\bullet )$ omitting the pull-back notation for $\mcal{C}$, when it is convenient.

Allowing the extra $Y$ turns out useful in our related work \cite{IwasaKai} with Ryomei Iwasa where $Y$ will be the projective space $\bbP ^r$.

%%%%%%%%%%%%%
%% Theorem %%
%\begin{definition}
%Let $((X,D),\mcal{C})$ be a pair $(X,D)$ and a family $\mcal{C}$ of constructible subsets of $X\setminus D$, and $Y$ an equi-dimensional scheme.
%We define the presheaf of cubical abelian group
%%
%\begin{equation*}
%\underline{z}^{iY}_{\mcal{C},X|D}(\bullet )
%\quad \text{ on }X_{\text{\'{e}t}}
% \end{equation*}
%%
%by $(U\xrightarrow{\phi }X)\mapsto \underline{z}^i_{\phi ^*\mcal{C}}(U\times Y|\phi ^*D\times Y,\bullet )$.
%When $Y=\Spec (k)$ or $\mcal{C}$ is trivial, we drop them respectively from the notation.

%\end{definition}
%% Theorem ends %%
%%%%%%%%%%%%%%%%%%

%%%%%%%%%%%%%%%%%%%%%%%%%%%%%%%
%%%%%%%%%%%%%%%%%%%%%%%%%%%%%%%
\subsubsection{}

The following is our main result in the full strength.

%%%%%%%%%%%%%
%% Theorem %%
\begin{theorem}\label{Th:full-statement}
Let $(X,D)$ be a pair of an equi-dimensional $k$-scheme over a field $k$ and an effective Cartier divisor, and $\mcal{C}$ an increasing family {\em (Definition \ref{Def:increasing-family})} of constructible subsets of $X\setminus D$.
Assume $X\setminus D$ is smooth.
Let $Y$ be an equi-dimensional $k$-scheme.
Then the inclusion map of the presheaves on $X_{\mrm{Nis}}$:
\begin{equation*}
\underline{z}^{i}_{\mcal{C}}(-\times Y\mid D_-\times Y,\bullet )
\to \underline{z}^{i}(-\times Y\mid D_-\times Y,\bullet )
 \end{equation*}
is a quasi-isomorphism in the Nisnevich topology.

\end{theorem}
%% Theorem ends %%
%%%%%%%%%%%%%%%%%%

Now we make some remarks before ending the section.

First, Theorem \ref{Th:full-statement} for a finite field $k$ follows from the case of infinite fields by the usual trace (norm) argument using extensions of $k$ having coprime degrees.
%See e.g.\ the first lines of \cite[Prop.4.7]{KrishnaParkAdditive} for a detailed discussion.
%A slight difference in our case is that 
We have to be careful in this case because
we are considering the stalks of presheaves. But it turns out that this does not pose a problem in view of Hensel's lemma saying that a Nisnevich neighborhood of a point of $X_{k'}$ (where $k\subset k'$ is a finite separable extension) can be refined by a Nisnevich neighborhood defined over $X$ up to the passage to connected components.

Second, Theorem \ref{Th:full-statement} implies the same assertion for the non-degenerate complexes $z^i$
because they are direct summands of the naive complexes $\underline{z}^i$ in a functorial way.
Of course we are primarily interested in the non-degenerate complexes.
Because of this, and to simplify the notation, in the text below we state results with underlines $\underline{z}^i$ but often omit them in the arguments.

Lastly, the following variant of Theorem \ref{Th:full-statement} implies the contravariant functoriality of additive higher Chow groups of smooth affine varieties. Their proofs are parallel, so we will mainly discuss the proof of Theorem \ref{Th:full-statement} in what follows and indicate necessary adaptation of the argument along the way as we need.

%%%%%%%%%%%%%
%% Theorem %%
\begin{theorem}\label{Th:additive-Chow}
Let $D\subset \bbA ^n_k$ be an effective Cartier divisor and $X$ a smooth affine $k$-scheme.
Suppose a family $\mcal{C}$ of constructible subsets of $X$ is given.

Then the inclusion
\begin{equation*}
\undl{z}^i_{\mcal{C}\times \bbA ^n}(X\times \bbA ^n\mid X\times D,\bullet )
\hookrightarrow 
\undl{z}^i(X\times \bbA ^n\mid X\times D,\bullet )
 \end{equation*}
is a quasi-isomorphism of complexes of abelian groups.

Consequently, additive higher Chow groups are contravariant in smooth affine varieties.
\end{theorem}
%% Theorem ends %%
%%%%%%%%%%%%%%%%%%

Here, $\mcal{C}\times \bbA ^n$ is the family $\{ C_{d-n}\times \bbA ^n \} _{d\in \bbZ }$ of
constructible subsets of
$X\times \bbA ^n$.

%%%%%%%%%%%%%%%%%%%%%%%%%%%%%%%
%%%%%%%%%%%%%%%%%%%%%%%%%%%%%%%
\section{Affine case}\label{Sec:Affine-case}

In this section, we prove the following special case of Theorem \ref{Th:full-statement}:

%%%%%%%%%%%%%
%% Theorem %%
\begin{theorem}\label{Th:affine-case}
Let $K$ be a discrete valuation field with valuation ring $\OF $,
and $\pi \in \OF$ any element in the maximal ideal.
Then for any sequence $\mcal{C}$ of constructible subsets of $\bbA ^n_{\OF}$ and any scheme $Y$ equi-dimensional over $\OF $,
the inclusion of complexes of abelian groups
\begin{equation*}
\underline{z}^i_{\mcal{C}}(\bbA ^n_{\OF} \times _{\OF }Y\mid V(\pi ),\bullet )
\subset 
\underline{z}^i(\bbA ^n_{\OF} \times _{\OF}Y\mid V(\pi ),\bullet )
 \end{equation*}
is a quasi-isomorphism.
Here, we denote by $V(\pi )$ the Cartier divisor defined by $\pi $.
\end{theorem}
%% Theorem ends %%
%%%%%%%%%%%%%%%%%%

We also record the following slightly more global version. It is useful in the context of additive higher Chow groups.

%%%%%%%%%%%%%
%% Theorem %%
\begin{theorem}\label{Th:affine-additive-case}
Let $k$ be a field and
$\pi (x_1,\dots ,x_n)\in k[x_1,\dots ,x_n]$ be a non-zero polynomial, which defines an effective Cartier divisor $D:=\{ \pi (\bm{x}) =0 \} $ on $\bbA ^n_k$.
Suppose a family $\mcal{C}$ of constructible subsets of $\bbA ^n_k\setminus D$ and
an equi-dimensional $k$-scheme $Y$ are given.

Then the inclusion of complexes of abelian groups
\begin{equation*}
\underline{z}^i_{\mcal{C}}(\bbA ^n_k \times Y\mid D\times Y,\bullet )
\subset 
\underline{z}^i(\bbA ^n_k \times Y\mid D\times Y,\bullet )
 \end{equation*}
is a quasi-isomorphism for every $i$.
\end{theorem}
%% Theorem ends %%
%%%%%%%%%%%%%%%%%%

In the proof of Theorem \ref{Th:full-statement} below, we will indicate some modifications needed for Theorem \ref{Th:affine-case}.

%%%%%%%%%%%%%%%%%%%%%%%%%%%%%%%
%%%%%%%%%%%%%%%%%%%%%%%%%%%%%%%
\paragraph{The homotopy operator}

We shall use the following homotopy operator.
Let $\bm{v}\in \bbA ^n({\OF})$ be  vector and $s\ge 0$ be an integer.
Consider the map of $\OF $-schemes
\begin{equation*}
\begin{array}{rccl}
\Phi _{\bm{v},s}\colon  &\bbA^n_{\OF}\times \bbA ^1 & \to &\bbA ^n_\OF  \\
&(\bm{x},t) &\mapsto &\bm{x}+\pi ^st \bm{v} .
 \end{array}
 \end{equation*}
When $t=0$, the map $\Phi _{\bm{v},s}(-,0)\colon $ $\bbA ^n_\OF \to \bbA ^n_\OF $ 
is the identity map and when $t=1$, the map $\Phi _{\bm{v},s}(-,1)$
is the translation by the vector $\pi ^s\bm{v}$.

For every $p\ge 0$, let us identify $\bbA ^1\times \square ^p $ and $\square ^{p+1}$ by identifying $\bbA ^1$ with the last factor of $\square ^{p+1}$.
Then we can ask if the map
\begin{equation*}
\Phi _{\bm{v},s}\times \id _Y\times \id _{\square ^p}\colon
\bbA ^n_\OF \times \bbA ^1\times Y\times \square ^p
\to \bbA ^n_\OF \times Y\times \square ^p
 \end{equation*}
pulls back a given element of $\underline{z}^i(\bbA ^n_\OF \times Y\mid V(\pi ),p) $ 
into $\underline{z}^i_{\mcal{C}}(\bbA ^n_\OF\times Y\mid V(\pi ),p+1) $.
Since the face condition is easily seen to be always satisfied,
the problem is the modulus condition and the proper intersection with the sets $C_d$.
We are going to show that it can be made true by choosing suitable $s$ and $\bm{v}$.

[For Theorem \ref{Th:affine-additive-case}, 
we choose $\bm{v}$ as an element of $\bbA ^n(k) $ and define the map $\Phi _{\bm{v},s}\colon \bbA ^n_k\times \bbA ^1\to \bbA ^n_k$ by the same formula $(\bm{x},t)\mapsto \bm{x}+t\pi (\bm{x})^s \bm{v}$
and consider pull-back by $\Phi _{\bm{v},s}\times \id _Y \times \id _{\square ^p}$:
$\bbA ^n_k\times \bbA ^1\times Y\times \square ^p $ $\to $ $\bbA ^n_k\times Y\times \square ^p $.]

We keep the notation throughout this section.
We omit the underline in $\underline{z}^i$ below to simplify the notation.

%%%%%%%%%%%%%%%%%%%%%%%%%%%%%%%
%%%%%%%%%%%%%%%%%%%%%%%%%%%%%%%
\subsection{Verifying the modulus condition}

Here we prove the following.

%%%%%%%%%%%%%
%% Theorem %%
\begin{proposition}\label{Prop:s(V)}
Let $V\subset \bbA ^n_\OF \times Y\times \square ^p$ be an irreducible closed subset satisfying the modulus condition with respect to the divisor $V(\pi )$.
Then there is an integer $s(V)\ge 0$ such that the closed set 
\begin{equation*}
\Phi \inv _{\bm{v},s}(V)\subset \bbA ^n_\OF \times Y\times \square ^{p+1}
 \end{equation*}
satisfies the modulus condition for every $s\ge s(V)$ and $\bm{v}\in \bbA ^n(\OF )$.
\end{proposition}
%% Theorem ends %%
%%%%%%%%%%%%%%%%%%

[This proposition is true for the case of Theorem \ref{Th:affine-additive-case} as it is.]

Since the problem is local on $Y$, we may assume $Y$ is affine.
Write $\Phi := \Phi _{\bm{v},s}$ and let
\begin{equation*}
\Phi \inv (V)^- \subset \bbA ^n_\OF \times Y\times (\bbP ^1) ^{p+1} \end{equation*}
be the closure of $\Phi \inv (V)$. For each index $i\in \{ 1,\dots ,p+1 \}$, denote by $\bbP ^1_i$ the $i$-th factor of $(\bbP ^1)^{p+1}$.
For a subset $I\subset \{ 1,\dots ,p+1 \} $, let us define
\begin{equation*}
U_I:= \prod _{i\in I}(\bbP ^1_i\setminus \{ 0\} )\times \prod _{i'\in I^c}(\bbP ^1_{i'}\setminus \{ \infty \} )
 \end{equation*}
which is an open set of $(\bbP ^1)^{p+1}$.
On the set $U_I$, the functions $1/t_i$ ($i\in I$) and $t_{i'}$ ($i'\in I^c$) are regular functions.
The divisor $F_{p+1}$ is defined by $\disp \prod _{i\in I}\frac{1}{t_i}$
in this region.

We define $U'_{I'}\subset (\bbP ^1)^p$ for a subset $I'\subset \{ 1,\dots ,p \} $ in a similar way.

Since the modulus condition is a local condition on the closed subset in question,
let us consider the modulus condition for 
\begin{equation*}
\Phi \inv (V)^-_I:=\Phi \inv (V)^- \times _{(\bbP ^1)^{p+1}}U_I
 \end{equation*}
for subsets $I\subset \{ 1,\dots ,p+1  \} $.

If $p+1\not\in I$, then the set $\Phi \inv (V)^-_I$ satisfies the modulus condition for every $s\ge 0$ and $\bm{v}$ because it is the pull-back of $\bar{V}\times _{(\bbP ^1)^p} U'_I$
(which satisfies the modulus condition)
by the map
\begin{equation*}
\Phi \times \id _{Y}\times \id _{(\bbP ^1)^p}\colon
\bbA ^n_\OF \times Y\times U_I
\to \bbA ^n_\OF \times Y\times U'_I .
 \end{equation*}

\newcommand{\Iminp}{{I' }}

Now suppose $p+1\in I$ and write $I':= I\setminus \{ p+1 \} $ which we often regard as a subset of $\{ 1,\dots ,p \} $.
Let $\mcal{I}$ be the ideal defining the integral closed subscheme 
\begin{equation*}
\bar{V}_\Iminp := \bar{V}\times _{(\bbP ^1)^p} U'_{\Iminp }  
\subset \bbA ^n_\OF \times Y\times U'_{\Iminp }
\end{equation*}
and let $\{ f_\lambda (\bm{x},\frac{1}{t_i} {}_{(i\in \Iminp )}, t_{i'}\ {}_{(i'\in I^c)}) \} _{\lambda \in \Lambda }$
be a finite set that generates $\mcal{I}$
(recall that we are assuming that $Y$ is affine).

Let us note the following trivial fact:
%%%%%%%%%%%%%
%% Theorem %%
\begin{lemma}\label{Lem:Modulus-homogeneous}
Let $A$ be an integral domain and $\Frac (A)$ be its fraction field.
Let $z,\pi \in A$ be two elements. Then the element ${z}/{\pi }$ of $\Frac (A)$ is integral over $A$ if and only if there is a homogeneous polynomial
\begin{equation*}
E(\alpha ,\beta )\in A[\alpha ,\beta ]
 \end{equation*}
which is monic as an element of $(A[\beta ])[\alpha ]$ such that
\begin{equation*}
E(z,\pi )=0 \quad \text{ in }A .
 \end{equation*}
\end{lemma}
%% Theorem ends %%
%%%%%%%%%%%%%%%%%%
(This is just the matter of clearing the denominators of an integral dependence relation
$\left( \frac{z}{\pi }\right) ^d+a_1\left( \frac{z}{\pi }\right)^{d-1}+\dots + a_d =0$ in $\Frac (A)$, $a_j\in A$.)

Now, since $\bar{V}_I$ satisfies the modulus condition,
Lemma \ref{Lem:Modulus-homogeneous} applied to $\bar{V}_I$
implies that we have a relation of the form
\begin{equation}\label{Eq:rel-on-V}
E_{I} [\bm{x},\frac{1}{t_i},t_{i'}](\prod _{i\in \Iminp }\frac{1}{t_i} ,\pi )
= \sum _{\lambda \in \Lambda } b_\lambda f_\lambda 
\end{equation}
where $E_{I} [\bm{x},\frac{1}{t_i},t_{i'}] (\alpha ,\beta )$
$\in \mcal{O}(\bbA ^n_\OF \times Y\times U'_\Iminp )[\alpha ,\beta ]$
is homogeneous in $\alpha ,\beta $ and monic in $\alpha $,
and $b_\lambda \in \mcal{O}(\bbA ^n_\OF \times Y\times U'_\Iminp )$.

%Substituting $\bm{x}+\pi ^st\bm{v}$ for $\bm{x}$ in the above formula gives
%%
%\begin{align*}
%E_I[\bm{x}+\pi ^st\bm{v},\frac{1}{t_i}, t_{i'} ]&(\prod _{i\in \Iminp }\frac{1}{t_i} ,\pi )
%\\
%&= \sum _{\lambda \in \Lambda } 
%b_\lambda (\bm{x}+\pi ^st\bm{v},\frac{1}{t_i}, t_{i'})
%f_\lambda (\bm{x}+\pi ^st\bm{v},\frac{1}{t_i}, t_{i'})
% \end{align*}
%%
%in the coordinate ring of $\bbA ^n_\OF \times \bbA ^1 \times Y\times U'_\Iminp $.

Meanwhile, let us observe the following.
Let $\deg _{\bm{x}}(f_\lambda ) $ be the total degree of $f_\lambda $
with respect to the variables $\bm{x}$. Then as a general fact about closure, %\cite[EGA]{EGA}, 
the regular functions
\begin{align*}
\phi _\lambda &(\bm{x},\frac{1}{t_i} {}_{(i\in I')}, t_{i'} {}_{(i'\in I^c)}, \frac{1}{t})
\\
&:= f_\lambda (\bm{x}+\pi ^st\bm{v},\frac{1}{t_i}, t_{i'})
\cdot \left( \frac{1}{t}\right) ^{\deg _{\bm{x}}(f_\lambda )}
 \end{align*}
on $\bbA ^n_{\OF} \times Y\times U_{I}$ are in the defining ideal for $\Phi \inv (V)^-_I$ 
because it is true on the locus where both $t$ and $\frac{1}{t}$ are regular.

Expanding the right hand side, it can be written as
\begin{align}\label{Eq:expand-f}
\phi _\lambda (\bm{x},\frac{1}{t_i},t_{i'},\frac{1}{t}) 
&= f_\lambda (\bm{x},\frac{1}{t_i}, t_{i'})
\cdot \left( \frac{1}{t}\right) ^{\deg _{\bm{x}}(f_\lambda )}\notag
\\
&\phan{=}+\pi ^s g_\lambda (\bm{x},\frac{1}{t_i}, t_{i'},\frac{1}{t}) 
 \end{align}
for some $g_\lambda \in \mcal{O}(\bbA ^n_\OF \times Y\times U_I)$.

Now, let $\deg (E_I)$ be the homogeneous degree of $E_I$ and set
\begin{equation*}
s(V):= \max (\{ \deg (E_I) \} _I \cup \{ \deg _{\bm{x}}(f_\lambda ) \} _{I,\lambda } ).
 \end{equation*}
This depends only on $V$ and not on $s$ or $\bm{v}$.
Suppose $s\ge s(V)$.
Then we may multiply \eqref{Eq:expand-f} by $\disp \left( \frac{1}{t}\right) ^{s-\deg _{\bm{x}}(f_\lambda ) } b_\lambda (\bm{x},\frac{1}{t_i},t_{i'})$ and sum them up with respect to $\lambda \in \Lambda $;
formula 
\eqref{Eq:rel-on-V}
then gives
\begin{align}\label{Eq:one-step-away}
E_I
&(\prod _{i\in I'}\frac{1}{t_i},\pi )
\cdot \left( \frac 1 t \right) ^s
+\pi ^s g(\bm{x},\frac{1}{t_i},t_{i'},\frac{1}{t})
\\
&= \sum _{\lambda \in \Lambda }
b_\lambda (\bm{x},\frac{1}{t_i},t_{i'})\cdot \left( \frac 1 t \right) ^{s- deg _{\bm{x}}(f_\lambda ) }\cdot
\varphi _\lambda (\bm{x},\frac{1}{t_i},t_{i'},\frac{1}{t}) \notag
 \end{align}
(with $g(\bm{x},\frac{1}{t_i},t_{i'},\frac{1}{t})\in \mcal{O}(\bbA ^n_\OF \times Y\times U_I)$)
whose right hand side is in the defining ideal for $\Phi \inv (V)^-_I$.

Write $E_I$ as the sum of its leading term and the other:
$E_I(\alpha ,\beta )=\alpha ^{\deg (E_I)}+ {E}'_I(\alpha ,\beta )$.
Then the left hand side of formula \eqref{Eq:one-step-away} multiplied by $\disp \left( \prod _{i\in I'}\frac{1}{t_i}\right) ^{s-\deg (E_I)}$ reads
\begin{equation*}
\left( \frac{1}{t} \prod _{i\in I }\frac{1}{t_i} \right) ^s
+ 
{E}'_I( \prod _{i\in I }\frac{1}{t_i} ,\pi  )\cdot \left( \frac{1}{t}  \right) ^s\left(  \prod _{i\in I '}\frac{1}{t_i} \right) ^{ s-\deg (E_I) }
+\pi ^s h
 \end{equation*}
with $h\in \mcal{O}(\bbA ^n_\OF \times Y\times U_I)$.
Clearly, this can be written as
$\tilde{E}_I(\prod _{i\in I}\frac{1}{t_i},\pi )$
for an $\tilde{E}_I[\bm{x},\frac{1}{t_i},t_{i'},\frac{1}{t}](\alpha ,\beta )\in \mcal{O}(\bbA ^n_\OF \times Y\times U_I )[\alpha ,\beta ]$
which is homogeneous in $\alpha ,\beta $ of degree $s$ and monic in $\alpha $.
By Lemma \ref{Lem:Modulus-homogeneous} again (applied to irreducible components of $\Phi \inv (V)^-_I$), we conclude that the closed subset $\Phi \inv (V)^-_I$ satisfies the modulus condition on this region. Since this is true for all $I$, Proposition \ref{Prop:s(V)} has been proven.

\ \ \

[As is clear from the proof, this choice of $s(V)$ works universally, in the sense that the conclusion of the lemma is true after any extension of DVRs $\OF \subset \mcal{O}'$, though we do not need this.]

%%%%%%%%%%%%%
%% Theorem %%
\begin{remark}\label{Rem:cheaper-way}
There is a way of proving Theorem \ref{Th:affine-case}
without proving Proposition \ref{Prop:s(V)}.
It consists of
pulling back a cycle by $\Phi _{\bm{v},1}$ with $\bm{v}$ sufficiently general (see the next subsection) and pushing it forward by the $d$-th power map $\square ^1\to \square ^1$ (in the last factor of $\square ^{p+1}$) with high enough powers $d$; see \cite[\S 4.2.3]{Miyazaki}.

Nonetheless, I chose to preserve Proposition \ref{Prop:s(V)} 
partly because my method has not appeared in the literature and
partly because the method has been applied in \cite{KrishnaParkCrys} for their ``fs-moving lemma.''

\end{remark}
%% Theorem ends %%
%%%%%%%%%%%%%%%%%%

%%%%%%%%%%%%%%%%%%%%%%%%%%%%%%%
%%%%%%%%%%%%%%%%%%%%%%%%%%%%%%%
\subsection{Verifying proper intersection}\label{Sec:proper-intersection}

We keep our notation: Namely, let $\OF $ be a discrete valuation ring and $\pi \in \OF $ an element in the maximal ideal.
Let $\mcal{C}=\{ C_d \} _d$ be a sequence of constructible subsets on $\bbA ^n_K$
and $Y$ be a scheme equi-dimensional over $\OF $, with fiber dimension
$\dim (Y/\OF )$.
Let $s\ge 0$ be an integer and the symbol $\bm{v}$ denote elements of $\bbA ^n(\OF )$.
We are interested in pull-back of cycles by the morphism
\begin{equation*}
\begin{array}{rccc}
\Phi _{\bm{v},s}\times \id _Y\times \id _{\square ^p}\colon
&\bbA ^n_\OF \times _\OF Y \times \square ^p \times \bbA ^1
&\to &\bbA ^n_\OF \times _\OF Y \times \square ^p 
\\
&(\bm{x},y,\bm{t},t)&\mapsto &(\bm{x}+\pi ^st\bm{v},y,\bm{t}) .
 \end{array}
 \end{equation*}
It is useful to allow more general vectors $\bm{v}$ when we discuss properties of a {\it general} $\bm{v}$. Namely,
for any ring homomorphism $\OF \to A  $ and a vector $\bm{v}\in \bbA ^n(A)$, this morphism makes sense as a morphism from
$\bbA ^n_A \times _\OF Y\times \square ^p \times \bbA ^1$.
We are primarily interested in the case where $A$ is the fraction field $K$ or its extension.

When we say some condition is true for a {\it general} $\bm{v}$, 
we mean that there is a dense open subset $U\subset \bbA ^n_K$ such that 
the condition is true for every vector of the form $\bm{v}\colon \Spec (\Omega )\to U $, with $\Omega $ a field.
Here, we concern only the generic fiber $\bbA ^n_K$ 
because the conditions we consider will depend only on what happens over the generic point of $\Spec (\OF)$.

%%%%%%%%%%%%%
%% Theorem %%
\begin{proposition}\label{Prop:proper-intersection-affine}
Let $s\ge 0$ be an integer
and 
$V\subset \underline{z}^i(\bbA ^n_\OF \times _\OF Y\mid V(\pi ), p)$ be a cycle. 
Then a general $\bm{v}\in \bbA ^n_K$ satisfies the following:
For every face $F $ of $\square ^p\times \bbA ^1\cong \square ^{p+1}$
and integer $d\in \bbZ $,
we have the inequality of dimension:
\begin{align*}
\dim &([C_d\times _\OF Y\times (F \setminus \{ t=0\} ) ]\cap |\Phi ^*_{\bm{v},s}(V)|) \\
\le &[\dim (C_d)+\dim (Y/\OF )+ \dim (F)] -i.
 \end{align*}
Consequently, if this $\bm{v}$ comes from an element of $\bbA ^n(\OF )$,
it follows that the cycle $\Phi _{\bm{v},s}^*(V)|_{t=1}$
is  in the subgroup
$\underline{z}^i_{\mcal{C}}(\bbA ^n_\OF \times _\OF Y\mid V(\pi ),p)$;
if moreover $V$ happens to be in the subgroup
$\underline{z}^i_{\mcal{C}}(\bbA ^n_\OF \times _\OF Y\mid V(\pi ),p)$,
then 
$\Phi _{\bm{v},s}^*(V)$ is in $\underline{z}^i_{\mcal{C}}(\bbA ^n_\OF \times _\OF Y\mid V(\pi ),p+1)$.

\end{proposition}
%% Theorem ends %%
%%%%%%%%%%%%%%%%%%

The following proof is essentially the same as that for a criterion given in \cite[Lem.1.1]{Bloch1986} which is often cited in this area. Here it seems to be more reader-friendly to give a direct proof than to make an obscure coordinate change so that we can apply the criterion.

Since there are only finitely many faces in a cube and data of (non-trivial) subsets $C_d$, we fix a face $F$ of $\square ^p\times \bbA ^1$
and an integer $d$.
Let $\Sigma $ be the subset of $\bbA ^n_K\times \bbA ^n_K\times (\bbA ^n_K\times Y_K\times (F\setminus \{ t=0\} ))$
consisting of points $(\bm{v}, \bm{x}', (\bm{x},y,\bm{t},t))$ satisfying the following four conditions:
\begin{itemize}
\item $\bm{x}'\in C_d$;
\item $(\bm{x},y,\bm{t})\in |V| $;
\item $\bm{x}=\bm{x}'+\pi ^st\bm{v}$.
 \end{itemize}
We want to consider $\Sigma $ because the set $[C_d\times _\OF Y\times (F \setminus \{ t=0\} ) ]\cap |\Phi ^*_{\bm{v},s}(V)|$ can be recovered as follows:
Look at the fiber $\Sigma _{\bm{v}}$ of the projection to the $\bm{v}$-factor $\Sigma \to \bbA ^n_K$ above the chosen point $\bm{v}$ and 
take the image by the projection $\Sigma _{\bm{v}}\to \bbA ^n_K\times Y_K\times (F\setminus \{ t=0\} )$ to the $(\bm{x}',y,\bm{t},t)$-factors;
then this is the set we are interested in.

Now, for every pair of
$\bm{x}'\in C_d $ and $(\bm{x},y,\bm{t},t)\in (|V|\times \bbA ^1)\times _{\square ^{p+1}} (F\setminus \{ t=0\} )$, there is a unique $\bm{v}\in \bbA ^n_K$ that makes the third condition true (namely $\bm{v}:= (\bm{x}-\bm{x}')/\pi ^s t$).
It follows that 
\begin{align*}
\dim (\Sigma ) 
&=\dim (C_d)+\dim ((|V|\times \bbA ^1)\times _{\square ^p+1} (F\setminus \{ t=0\} ))
\\
&\le \dim (C_d)+\dim (\bbA ^n_K\times Y_K \times F) -i
\quad \text{ [the face condition]}.
\end{align*}
Therefore, the projection to the $\bm{v}$-factor
$\Sigma \to \bbA ^n_K$ has generic fiber with dimension
$\le \dim (C_d)+\dim (Y_K\times F) -i$, cf.\ Proposition \ref{Cor:Chevalley-dim} (Chevalley's theorem).
This complets the proof of Proposition \ref{Prop:proper-intersection-affine}.

\ \ \

[For Theorem \ref{Th:affine-additive-case},
one needs the inequality 
\begin{align*}
\dim &([C_d\times  Y\times (F \setminus \{ t=0\} ) ]\cap |\Phi ^*_{\bm{v},s}(V)|) \\
\le &[\dim (C_d)+\dim (Y )+ \dim (F)] -i.
 \end{align*}
The above proof works if one defines $\Sigma \subset \bbA ^n_k\times \bbA ^n_k\times (\bbA ^n_k\times Y\times (F\setminus \{ t=0 \}))$ 
by the same formula.]

%%%%%%%%%%%%%%%%%%%%%%%%%%%%%%%
%%%%%%%%%%%%%%%%%%%%%%%%%%%%%%%
\subsection{Proof of the affine case}\label{Sec:proof-of-affine-case}

Now we can prove that the inclusion 
\begin{equation*}
\underline{z}^i_{\mcal{C}}(\bbA ^n_\OF \times _\OF Y\mid V(\pi ),\bullet )\to 
\underline{z}^i(\bbA ^n_\OF \times _\OF Y\mid V(\pi ),\bullet )
 \end{equation*}
induces isomorphisms on the homology groups.

For the surjectivity of the map on the homology, choose any homology class in the target and a cycle $V\in \underline{z}^i(\bbA ^n_\OF \times _\OF Y\mid V(\pi ),p )$ representing it. 
Take a large enough $s\ge 0$ using Proposition \ref{Prop:s(V)}
and then a vector $\bm{v}\in \bbA ^n_\OF $ whose generic fiber is sufficiently general using Proposition \ref{Prop:proper-intersection-affine}.
We conclude that the cycle $\Phi _{\bm{v},s}^*(V)|_{t=1}$ is in the subcomplex
$\underline{z}^i_{\mcal{C}}(\bbA ^n_\OF \times _\OF Y\mid V(\pi ),p )$
and the obvious compatibility tells us
\begin{equation*}
\Phi ^*_{\bm{v},s}(V)|_{t=1}-V
= (-1)^p\partial (\Phi ^*_{\bm{v},s}(V) )
 \end{equation*}
where $\partial $ denotes the differential in the cycle complex.
It follows that the map on the homology is surjective.

The proof of the injectivity is similar: We take any cycle $W$ representing a homology class in the source, suppose it becomes a boundary $W=\partial V$ in the target and choose $s\ge 0$ and $\bm{v}$ such that $\Phi _{\bm{v},s}^* $ works for $V$. The point is that the inclusion in question is a {\it weak homotopy equivalence}, meaning that even if we might not have a uniform homotopy inverse, we do have a homotopy inverse on every finitely generated subcomplex.

This completes the proof of Theorem \ref{Th:affine-case}.

[The corresponding step for the proof of Theorem \ref{Th:affine-additive-case} is the same.]

%%%%%%%%%%%%%%%%%%%%%%%%%%%%%%%
%%%%%%%%%%%%%%%%%%%%%%%%%%%%%%%
\section{Linear projection method}\label{Sec:linear}

We collect some general facts on the technique of ``moving by linear projections''
$\bbA^N\twoheadrightarrow \bbA^n $ ($N\ge n$).

Let us write $\bbA^N=\Spec (\bbZ [x_1,\dots ,x_N])$
and $\bbA ^n=\Spec (\bbZ [y_1,\dots ,y_n])$.

%%%%%%%%%%%%%%%%%%%%%%%%%%%%%%%
%%%%%%%%%%%%%%%%%%%%%%%%%%%%%%%
\subsection{Notation}

Suppose we are given $n$ linear polynomials
$\phi _1,\dots ,\phi _n $ in $N$ variables
(i.e.\ an $N\times n$ matrix) with coefficients in some ring.
They define a map of schemes $\phi \colon \bbA^N\to \bbA ^n$ over the ring,
defined by $\bm{x}=(x_1,\dots ,x_N)\mapsto (\phi _1(\bm{x}),\dots ,\phi _n(\bm{x}))$.
The surjective linear maps form a dense open subset $M_{N\times n}^*$ of the space of all $N\times n$ matrices $M_{N\times n}$.

The following projective terminology is often useful as well.
Let $\bbA^N\subset \bbP ^N$ be the usual open immersion.
Let $X_0,\dots ,X_N$ be the homogeneous coordinates so that $x_i=X_i/X_0$ ($i=1,\dots ,N$).
We may think of the polynomials $\phi _i(x_1,\dots ,x_N)$ as the dehomogenization of the polynomials $\phi _i(X_1,\dots ,X_N)$.

Define a morphism $\Phi $ by
\begin{equation*}
\begin{array}{ccl}
\Phi \colon \bbP^N \setminus V((X_0,\phi _1(\bm{X}),\dots ,\phi _n(\bm{X})))
&\to &\bbP ^n  \\{}
[X_0:\dots :X_N]&\mapsto & {} 
[X_0: \phi _1(\bm{X}):\dots :\phi _n(\bm{X})]
\end{array}
 \end{equation*}
which is affine. %\cite[??]{EGA}.
Let us write $V(\phi ):=V((X_0,\phi _1(\bm{X}),\dots ,\phi _n(\bm{X})))$.
It can be checked that the following diagram:
\begin{equation*}
\xymatrix@R=10pt@C=50pt{
\bbP^N \setminus V(\phi )
\ar[r]^{\Phi }
& \bbP ^n
\\%
\bbA ^N\ar@{}|{\cup }[u]
\ar[r]_{\phi }
&\bbA ^n\ar@{}|{\cup }[u]
\ar@{}|{\square }[ul]
}
 \end{equation*}
is cartesian.

If $\phi $ belongs to $M^*_{N\times n}$, the morphism $\Phi $ is an affine bundle with fiber $\bbA ^{N-n}$.

%%%%%%%%%%%%%%%%%%%%%%%%%%%%%%%
%%%%%%%%%%%%%%%%%%%%%%%%%%%%%%%
\subsection{Over a base field}

From now on, we work over a base field $k$, and $\phi $ will always be taken from $M^*_{N\times n}$.

When we say some property holds for a general $\phi \in M^*_{N\times n}$, it means that there is a dense open subset of $M^*_{N\times n}$ such that if a field-valued point $\Spec (\Omega )\to M^*_{N\times n}$ maps into the subset, then the property holds for the corresponding map $\phi \colon \bbA^N_{\Omega }\to \bbA ^n_{\Omega }$.
Since the properties we consider are such that the points satisfying the property form a constructible subset of $M^*_{N\times n}$, 
this will always be equivalent to the condition that the property holds for 
the generic point of $M^*_{N\times n}$.
Also, since $M^*_{N\times n}$ is of finite type over the field $k$, it is equivalent to its truth for the closed points of a dense open subset of $M^*_{N\times n}$.
If $k$ is infinite, it can even be checked on $k$-rational points alone.

For a field extension $k\subset \Omega $, a matrix $\phi \in M^*_{N\times n}(\Omega )$ and a point $\bbA ^N(\Omega )$,
the fiber of $\phi \colon \bbA ^N_\Omega \to \bbA ^n_\Omega $ containing $x$ is $x + \ker (\phi )$.

Let $X\subset \bbA ^N$ be a closed subset and $\bar{X} $ be its closure in $\bbP^N$.
If $\bar{X}$ does not meet $V(\phi )$, then $\Phi _{| \bar{X}}$ is both affine and projective, hence finite (regardless of the scheme structure of $\bar{X}$).
It follows that $\phi _{|X}$ is finite and $\dim (X)\le n$ in this case.
Conversely, if $X$ has dimension $\le n$, this becomes the case for a general $\phi \in M^*_{N\times n}$.

%%%%%%%%%%%%%%%%%%%%%%%%%%%%%%%
%%%%%%%%%%%%%%%%%%%%%%%%%%%%%%%
\subsubsection{Avoiding a closed subset}

We keep the notation above.

Let $x\in \bbA^N$ be a point and $Z\subset \bbA ^N $ be a closed subset of dimension $<n$ not containing $x$.
Then for a general $\phi \in M^*_{N\times n}$,
we have $\phi (x)\not\in \phi (Z)$.

%%%%%%%%%%%%%
%% Theorem %%
\begin{proof}
Define a closed subset $\Sigma \subset Z\times M^*_{N\times n}$ by
\begin{equation*}
\Sigma := \left\{ (z,\phi )  \mid \phi (z)=\phi (x) \right\} .
 \end{equation*}
The fiber of a given $z\in Z$ along the projection $\Sigma \to Z$ is
\begin{equation*}
\left\{  \phi \in M^*_{N\times n} \mid  x\in z+\ker (\phi ) \right\}
 \end{equation*}
(i.e.\ the vector $\overrightarrow{xz}$ connecting $x$ and $z$ is in the kernel of $\phi $).
This set has codimension $n$ in $M^*_{N\times n}$ because $\overrightarrow{xz}\neq \bm{0}$ by the assumption $x\not\in Z$.
It follows that $\dim (\Sigma )=\dim (M^*_{N\times n})+\dim (Z)-n <\dim (M^*_{N\times n})$.
Therefore the projection $\Sigma \to M^*_{N\times n}$ cannot be surjective.
It follows that for the generic point $\phi \in M^*_{N\times n}$, there is no $z\in Z$ satisfying $\phi (z)=\phi (x)$.
\end{proof}
%% Theorem ends %%
%%%%%%%%%%%%%%%%%%

%%%%%%%%%%%%%%%%%%%%%%%%%%%%%%%
%%%%%%%%%%%%%%%%%%%%%%%%%%%%%%%
\subsubsection{Smoothness}\label{Sec:smoothness}

We need to recall Veronese embeddings.
Let $d\ge 1$ be an integer.
Let $V\subset k[x_1,\dots ,x_N]$ be the linear subspace of polynomials of degree $\le d$.
We will denote the affine space $\Spec (k[V])$ also by $V^*$.
The Veronese embedding of degree $d$ is the closed embedding
\begin{equation*}
\bbA ^N\hookrightarrow V^*
 \end{equation*}
corresponding to the surjective ring map
$k[V]\twoheadrightarrow k[x_1,\dots ,x_N]$ induced by the inclusion of $V$ into the right hand side.
Occasionally we consider the corresponding projective counterpart
\begin{equation*}
\bbP ^N\hookrightarrow \bbP (k\oplus V^*).
 \end{equation*}
The Veronese re-embedding of a given closed embedding $X\hookrightarrow \bbA ^N$ refers to the composition with the $d$-fold Veronese embedding
$X\rightarrow \bbA^N\hookrightarrow V^*$.
When we say some property holds for a general $\phi \in M^*_{N\times n}$ after the Veronese re-embedding of some degree, we mean that it is true after composing $X\hookrightarrow \bbA ^N\hookrightarrow V^*$ and writing $\bbA ^N $ anew for $V^*$. 

\ \ \

Let $X\subset \bbA ^N$ be a closed subscheme of dimension $\ge n$ everywhere, and $x\in X$ be a point.
Then the map $\phi _{|X}$ is smooth along $\phi ^{-1}\phi (x)\cap X_{\mrm{sm}}$ for a general $\phi \in M^*_{N\times n}$
after any Veronese re-embedding of degree $\ge 2$.
(If $\dim (X)=n$, we do not need a Veronese re-embedding.)

%%%%%%%%%%%%%
%% Theorem %%
%\begin{proof}
This is a version of Bertini smoothness theorem; see \cite[XI, 2.1]{SGA4}.
%\end{proof}
%% Theorem ends %%
%%%%%%%%%%%%%%%%%%

%%%%%%%%%%%%%%%%%%%%%%%%%%%%%%%
%%%%%%%%%%%%%%%%%%%%%%%%%%%%%%%
\subsubsection{Birationality}\label{Sec:birationality}

Let $X\subset \bbA^N$ be a closed subscheme of dimension $\le n$ and $x\in X_{\mrm{sm}}$ be a point such that its closure has dimension $< n$.
Then the extension of residue fields $k(\phi (x))\to k(x)$ is trivial for a general $\phi \in M^*_{N\times n}$.

%%% Proof begins %%%
\begin{proof}
We prove the separability and the pure inseparability of the field extension separately.

First, for a general $\phi $, consider the composition of $\phi _{|X}$ and the projection to the first $\dim (X)$ factors: 
\begin{equation*}
\bar{\phi } \colon X\xrightarrow{\phi _{|X}} \bbA^n\to \bbA ^{\dim (X)} 
\end{equation*}
which is \'{e}tale at $x$ by \S \ref{Sec:smoothness}.
In particular, the extension $k(\bar{\phi}(x))\to k(x)$ is finite and separable. Therefore, so is $k(\phi (x))\to k(x)$.

Next, let $k\to k^{\mrm{perf}}$ be the perfect closure of $k$. Since it is a purely inseparable extension, there is a unique point $\tilde{x}$ on $\bbA^N_{k^{\mrm{perf}}}$ above $x$.
Since $k^{\mrm{perf}}$ is a perfect field, the residue field $k^{\mrm{perf}}(\tilde{x})$
is separably generated over $k^{\mrm{perf}}$.

Put $Z:= \overline{\{ x\} } $.
Refreshing the notation, let $\bar{\phi }$ denote the composition 
\begin{equation*}
\bar{ \phi }\colon Z_{k^{\mrm{perf}}}\xrightarrow{\phi }\bbA ^n_{k^{\mrm{perf}}}\to \bbA ^{\dim (Z)}_{k^{\mrm{perf}}} 
 \end{equation*}
where the second map is the projection to the first $\dim (Z)$ factors.
By the first half of this proof, the extension $k^{\mrm{perf}}(\bar{\phi }(\tilde{x}))\to k^{\mrm{perf}}(\tilde{x})$ finite and separable.
Note that the assumption $\dim (Z)<n $ allows us to consider the function $y_{\dim (Z)+1}$. Its pull-back $\phi ^*(y_{\dim (Z)+1})$ has the form
$c_1x_1+\dots +c_Nx_N$ with $c_i\in k^{\mrm{perf}}$ (or in some extension thereof).
Since the functions $x_1,\dots ,x_N$ generate the field $k^{\mrm{perf}}(\tilde{x})$ over $k^{\mrm{perf}}(\bar{\phi }(\tilde{x}))$, by (the proof of) Primitive Element Theorem, a general $\phi \in M^*_{N\times n}$ makes the equality 
\begin{equation*}
k^{\mrm{perf}}(\tilde{x})=k^{\mrm{perf}}(\bar{\phi }(\tilde{x})) (\phi ^*(y_{\dim (Z)+1}))  \end{equation*}
true.
In particular one has $k^{\mrm{perf}}(\tilde{x})=k^{\mrm{perf}}(\phi (\tilde{x}))$

Now in the following diagram:
\begin{equation*}
\xymatrix{
\Spec (k^{\mrm{perf}}(\tilde{x}))
\ar[d]\ar[r]^{(1)}
&\Spec (k^{\mrm{perf}}(\phi (\tilde{x}) )) 
\ar[d]^{(2)}
\\%
\Spec (k(x))
\ar[r]_{(3)}
& \Spec (k(\phi (x)))
}
 \end{equation*}
the arrow $(1)$ is an isomorphism and $(2)$ is purely inseparable. It follows that the arrow $(3)$ is purely inseparable as well.

This completes the proof.

\end{proof}
%%% Proof ends %%%

%%%%%%%%%%%%%%%%%%%%%%%%%%%%%%%
%%%%%%%%%%%%%%%%%%%%%%%%%%%%%%%
\subsubsection{Chow's moving lemma}\label{Sec:Chow's-lemma}

Let $X\subset \bbA^N$ be a closed subscheme of pure dimension $n$ and
let $W\subset X$ be a constructible subset.

Suppose we are given an equi-dimensional scheme $Y$, an integer $n\ge 0$ and
a closed subset $V$ of $X\times Y\times \square ^n$.

For a map $\phi \in M^*_{N\times n}$ such that $\phi _{|X}$ is finite surjective,
define a new closed subset $\phi ^+ (V)$ as the closure of $( \phi _{|X}\times \id _{Y\times \square ^n} )^{-1}(\phi _{|X}\times \id _{Y\times \square ^n}) (V)\setminus V$ in $X\times Y\times \square ^n$.

Now we claim that the following inequality holds:
\begin{align*}
\dim (\phi ^+(V)\cap (W\times Y))
\le \max \{ \dim (V)+\dim (W)-\dim (X)&, 
\\%
\dim (V\cap (W\times Y))-1 & \}
 \end{align*}
for a general $\phi \in M^*_{N\times n}$, if $X$ is smooth at the images of the generic points of $V\cap (W\times Y)$.

%%% Proof begins %%%
\begin{proof}
Clearly, we may assume $n=0 $ because $Y\times \square ^n$ is equi-dimensional anyway.

Denote by $\ram (\phi _{|X}) \subset X$ the complement of the \'{e}tale locus of $\phi _{|X}$.

Put $A:= V\cap (W\times Y)\cap \ram (\phi _{|X}\times \id _Y )$.
Observe the obvious inclusion
\begin{equation*}
\phi ^+(V)\cap (W\times Y)\quad \subset \quad 
A
\cup [\phi ^+(V)\cap (W\times Y) \setminus A]
 \end{equation*}
We are going to show
\begin{itemize}
\item $\dim (A)\le \dim (V\cap (W\times Y)) -1$ and

\item $\dim (\phi ^+(V)\cap (W\times Y) \setminus A)\le \dim (V)+\dim (W)-\dim (X)$
 \end{itemize}
for a general $\phi \in M^*_{N\times n}$.

First, since $X$ is assumed to be smooth at the images of the generic points of the set $V\cap (W\times Y)$, 
a general $\phi _{|X}$ is \'{e}tale at these points.
This shows the first inequality.

Next, consider the following constructible set:
\begin{align*}
\ell &(V,W):= 
\\%
&\left\{ (\ell , y) \middle| \begin{array}{l}
\ell \text{ is a line in }\bbA ^N, y\in Y, \\%
\exists (x_1,x_2)\in (X\times W )\setminus \Delta _X  \text{ such that }
x_1, x_2\in \ell \text{ and }(x_1,y)\in V
\end{array}
 \right\} .
 \end{align*}
It is easy to see the inequality $\dim (\ell (V,W))\le \dim (V)+\dim (W)$.
We have the ``direction'' map
\begin{equation*}
\delta \colon \ell (V,W)\to \bbP _\infty ^{N-n}.
 \end{equation*}
Two points $x_1,x_2\in X$ are mapped by $\phi $ to the same point if and only if either $x_1=x_2$ or $x_1\neq x_2$ and $\delta (\overrightarrow{x_1x_2})\in V(\phi )$.

Look at the following diagram of constructible sets:
\begin{equation*}
\begin{array}{l}
X\supset \phi ^+(V)\cap (W\times Y) \setminus A \\%
\\%
\uparrow \pr _1 \\%
\\%
\disp \Sigma _1:= \left\{ (x,\ell )\in W\times \ell (V,W) \ \middle| \ x\in \ell \text{ and } \delta (\ell )\in V(\phi )  \right\}
\\%
\\%
\downarrow \pr _2
\\%
\\%
\delta ^{-1}(V(\phi ))\subset \ell (V,W)\xrightarrow{\delta }\bbP _\infty ^{N-1}.
 \end{array}
 \end{equation*}
Now, since a general $\phi _{|X}$ is finite, the projection $\pr _2$ has finite fibers.
It follows that $\dim (\Sigma _1)= \dim (\delta ^{-1}(V(\phi )))$.

For each integer $j\ge 0$, write $S_j$ for the set of points of $\bbP^{N-1}_\infty $ at which the fiber of $\delta $ has dimension $\ge j$. Then if $V(\phi )$ meet the subsets $S_j$ properly in $\bbP _\infty ^{N-1}$, the set $\delta ^{-1}(V(\phi ))$ has dimension $\ell (V,W)-n $ which is at most $\dim (V)+\dim (W)-n$.
This proper intersection is true for a general $\phi \in M^*_{N\times n}$ because it is about intersections of a linear subspace and subsets of $\bbP ^{N-1}_\infty $.
We obtain
\begin{equation*}
\dim (\Sigma _1)\le \dim (V)+\dim (W)-\dim (X).
 \end{equation*}

It remains to show that the image of $\pr _1$ contains $\phi ^+(V)\cap (W\times Y)\setminus A$, because then we will obtain
\begin{equation*}
\dim (\phi ^+(V)\cap (W\times Y)\setminus A)\le \dim (\Sigma _1 )
\end{equation*}
which will finish the proof.

Let us observe
that since $V$ is closed and $\phi _{|X}$ is finite (hence a closed map), the subset $\phi ^+(V)$ is alternatively obtained by the following steps:
consider the subset 
\begin{equation*}
\Sigma _2=\{ (z_1,z_2)\in ((X\times Y)\setminus V )\times V \mid \phi (z_1)=\phi (z_2) \} ,
\end{equation*}
take its closure $\overline{\Sigma _2}$ in $(X\times Y)^2$ and project it by the first projection $\overline{\Sigma _2}\hookrightarrow (X\times Y)^2\to X\times Y$.

%%%%%%%%%%%%%
%% Theorem %%
\begin{lemma}[{adaptation of \cite[Prop.4]{Roberts1971}}]\label{Lem:Roberts}
Keep the notation just above.
Let $z=(x,y)\in X\times Y$ be a point.
If the point $(z,z)\in (X\times Y)^2$  lies in $\overline{\Sigma _2}$,
then there is a line $\ell \subset T_{X,x}$ with the property that $\phi (\ell )\in V(\phi )$.
\end{lemma}
%% Theorem ends %%
%%%%%%%%%%%%%%%%%%
Let us see how Lemma \ref{Lem:Roberts} implies the inclusion
$\phi ^+(V)\cap (W\times Y)\setminus A\subset \image (\pr _1)$.
Suppose we are given a point $z_1=(x_1,y)\in W\times Y$ which is in $\phi ^+(V)$.
By its definition, there is a point $z_2=(x_2,y)\in V$ with $\phi (x_1)=\phi (x_2)$. 
If $z_1\neq z_2$, then $(z_1,\overrightarrow{z_1z_2})$ gives a point whose first projection is $z_1$.
If $z_1=z_2$, then Lemma \ref{Lem:Roberts} implies that there is a line $\ell \subset T_{X,x_1}$ with $\delta (\ell )\in V(\phi )$, which implies that $\phi _{|X}$ is not \'{e}tale at this point. It follows that the point $z_1=z_2$ is in $V\cap (W\times Y)\cap \ram (\phi _{|X}\times \id _Y)=A$.

\ \ \

Finally, we prove Lemma \ref{Lem:Roberts}.
By \cite[II 7.1.7]{EGA}, there is a discrete valuation ring $R$ and a morphism $\Spec (R)\to (X\times Y)^2$ such that the generic point maps into $\Sigma _2$ and the closed point maps to $(z,z)$.
By extending $R$, we may assume that $R$ is of the form $\Omega [[t]]$ where $\Omega $ is a field.

Let us denote the first and the second components of the map by
\begin{equation*}
z_i(t)=(\bm{x}_i(t),y_i(t))\colon \Spec (\Omega [[t]])\to X\times Y,\qquad  i=1,2; 
\end{equation*}
then we have in particular $z_1(0)=z_2(0)=z$ in $(X\times Y)(\Omega )$.

Denote the coordinates of $\bm{x}_i(t)$ in $\bbA ^N $ by
\begin{align*}
\bm{x}_i(t)&=(x_{i,1}(t),\dots ,x_{i,N}(t)), \quad \text{each of which are power series:}\\%
x_{i,j}(t)&=\sum _{h\ge 0}x_{i,j}^{(h)}t^h \in \Omega [[t]].
 \end{align*}
Write $\bm{x}_i^{(h)}:=(x_{i,1}^{(h)},\dots ,x_{i,N}^{(h)})\in \Omega ^N$.

Since the generic point of $\Spec (\Omega [[t]])$ maps into $\Sigma _2$, we have
$y_1(t)=y_2(t)$ in $Y(\Omega [[t]])$ and
\begin{equation*}
\bm{x}_1(t)\neq \bm{x}_2(t),\qquad \delta ({\bm{x}_1(t)-\bm{x}_2(t)}) \in V(\phi )_{\Omega ((t))}.
 \end{equation*}
In particular $\delta ({\bm{x}_1^{(h)}-\bm{x}_2^{(h)}})\in V(\phi )_{\Omega }$ for all $h$ such that the vector ${\bm{x}_1^{(h)}-\bm{x}_2^{(h)}}$ is non-zero (otherwise the left hand side does not make sense).

Let $h_0\ge 1$ be the integer such that the vectors $\bm{x}_1(t)$ and $\bm{x}_2(t)$ satisfy
$\bm{x}_1(t)\equiv \bm{x}_2(t)$ modulo $t^{h_0}$ and $\bm{x}_1(t)\not\equiv \bm{x}_2(t)$ modulo $t^{h_0+1}$.
We claim that the first non-vanishing vector ${\bm{x}_1^{(h_0)}-\bm{x}_2^{(h_0)}}\in \Omega ^N$ belongs to the subspace $T_{X,x}\otimes _{k(x)}\Omega $.

Let $f(x_1,\dots ,x_N)\in \mcal{O}(\bbA ^N)$ be any function belonging to the defining ideal for $X$.
We have $f(\bm{x}_i(t))=0$ for $i=1,2$ so that
\begin{equation}\label{Eq:equation-for-X}
f(\bm{x}_1(t))-f(\bm{x}_2(t)) =0 \quad \text{ in }\Omega [[t]].
 \end{equation}
Denote the coordinates of the point $\bm{x}_{1}(0)=\bm{x}_{2}(0)=x$ by $(x_1(0),\dots ,x_N(0))$.
Note that the functions $x_{i,j}(t)-x_j(0)$ are divisible by $t$ in $\Omega [[t]]$.
Expand the function $f$ with center at $x$, and isolate the degree $1$ terms:
\begin{equation*}
f=\left( \sum _{j=1}^N \frac{\partial f}{\partial x_j}(x)\cdot (x_j-x_j(0)) 
\right)
+f_{\ge 2} 
 \end{equation*}
where $f_{\ge 2}$ is the sum of degree $\ge 2$ terms 
with respect to $\{ x_j-x_j(0) \} _j$.
By the congruence $x_{i,j}(t)-x_j(0)\equiv 0 $ modulo $t$ and $x_{1,j}(t)\equiv x_{2,j}(t) $ 
modulo $t^{h_0}$, the power series
$f_{\ge 2}(\bm{x}_1(t))-f_{\ge 2}(\bm{x}_2(t))$ is divisible by $t^{h_0+1}$.
It follows that
formula \eqref{Eq:equation-for-X} modulo $t^{h_0+1}$ reads:
\begin{equation*}
\sum _{j=1}^N \frac{\partial f}{\partial x_j}(x)\cdot (x_{1,j}^{(h_0)}-x_{2,j}^{(h_0)}) t^{h_0} \equiv 0 \mod t^{h_0+1},
 \end{equation*}
which precisely means $({d}f_x) (\bm{x}_1^{(h_0)}-\bm{x}_2^{(h_0)})=0$ in $\Omega $.
Since this holds for all functions $f$ which vanish on $X$, the vector ${\bm{x}_{1}^{(h_0)}-\bm{x}_2^{(h_0)}}$ is tangent to $X$.
\end{proof}
%%% Proof ends %%%

%%%%%%%%%%%%%%%%%%%%%%%%%%%%%%%
%%%%%%%%%%%%%%%%%%%%%%%%%%%%%%%
\section{Noether's normalization lemma over a Dedekind base}\label{Sec:Noether}

The following is a variant of \cite[Th.\ 10.2.2]{LevineChow}
where he considers $X$ smooth over an arbitrary noetherian scheme $B$.
In the next section, it will be applied to the situation where 
we are given a principal effective Cartier divisor $D\subset X$ defined by an equation $\pi \colon X\to \bbA ^1=:B$.

%%%%%%%%%%%%%
%% Theorem %%
\begin{theorem}\label{Th:Noether}
Let $B$ be the spectrum of a Dedekind domain and $p\colon X\to B$ be an equi-dimensional morphism (i.e., it is a morphism of finite type where every irreducible component of $X$ dominates $B$ and all fibers have the same dimension).
Let $n$ be its fiber dimension.

Suppose we are given a point $x\in X$ and let $b=p(x)\in B$ be its image.
Assume $k(b) $ is an infinite field.
Then there exists Nisnevich neighborhoods
$(X',x)\to (X,x)$ and $(B',b)\to (B,b)$, a commutative diagram
\begin{equation*}
\xymatrix{
X'\ar[d]_{p'}\ar[r]
& X \ar[d]^{p}
\\%
B'\ar[r]
& B
}
 \end{equation*}
and a closed immersion $X'\hookrightarrow \bbA ^N_{B'}$ for some $N\ge 0$
such that if we denote by $\bar{X}'$ its closure in $\bbP ^N_{B'}$,
then $X'_b$ is dense in $(\bar{X}')_b$.
 
\end{theorem}
%% Theorem ends %%
%%%%%%%%%%%%%%%%%%

Consequently, after possibly replacing $B$ by a dense open subset, we can find a matrix $\phi \in M^*_{N\times n}(B)$ such that
$\phi _{|X'}\colon X'\to \bbA ^N_{B'}$ is finite and surjective.
%Hence Theorem \ref{} is a Nisnevich local version of Noether's Normalization Lemma over a Dedekind base.

\ \ \

In the proof below, given a map of schemes $X\to B$ and a property $\mathbf{P}$ of schemes (or of schemes and its subsets, when some subsets of $X$ are also given), we will say $\mbf{P}$ is true {\it fiberwise} over $B$ to mean that the property $\mbf{P}$ holds for every fiber of the morphism.

Many of the properties appearing below are such that if the property holds for the fiber over a point $b\in B$, then it holds for all fibers over the points in an open neighborhood of $b$.
In this case, if we arrange that the property be true for the fiber over a given $b$ and we replace $B$ by such an open subset $B_1$, then the property becomes true fiberwise over $B_1$.

%%% Proof begins %%%
\begin{proof}
If $n=0$, this follows from a version of Hensel's lemma. We proceed by induction on $n$, and let us assume $n\ge 1$ in what follows.

Embed $X$ into a projective space $\bbP ^N_B$ as a locally closed subscheme, and take its (reduced) closure $\bar{X}$.
Since $B$ is Dedekind, $\bar{X}$ is still equi-dimensional over $B$ because every injective ring homomorphism from a Dedekind domain into an integral domain is flat.

[This is the only point where we use the assumption that $B$ is Dedekind, except that it is of course built in the induction step as well. Probably there are many special cases where this equi-dimensionality can be verified for a more general base $B$.]

Note that $X$ is not necessarily {\it fiberwise} dense in $\bar{X}$ over $B$ because there can be an irreducible component of $\bar{X}_b$ contained in $\bar{X}\setminus X$, which is the main difficulty we are going to address.

Possibly after shrinking $B$ to an open neighborhood $B_1$ of $b$, we can choose a tuple of sections $\psi = (\psi _1,\dots ,\psi _n)\in \Gamma (\bbP ^N_{B_1}, \mcal{O}(1))^n$,
which forms part of a basis for $\Gamma (\bbP ^N_{B_1}, \mcal{O}(1))$,
such that their common zero-locus $V(\psi )$ does not meet $x$, and meets $\bar{X}_1:=X\times _{B}B_1 $ fiberwise properly over $B_1$, so that their intersection is finite and surjective over $B_1$.
Put $X_2:= (X\times _B B_1) \setminus V(\psi )$, which is an open neighborhood of $x$.

Let $\tilde{\bbP } ^N_{B_1}\xrightarrow{\pi } \bbP^N_{B_1}$ be the blow-up along $V(\psi )$.
Denote by $\tilde{X}_1$ the strict transform of $\bar{X}_1$.
It actually turns out that $\tilde{X}_1=\pi ^{-1}(\bar{X}_1)$, but we will only need the obvious inclusion $\tilde{X}_1\subset \pi ^{-1}(\bar{X}_1)$.
The rational map $\bbP ^N_{B_1}\dashrightarrow \bbP^{n-1}_{B_1}$: $[X_0:\dots :X_N]\mapsto [\psi _1(\bm{X}):\dots :\psi _n(\bm{X})]$
gives a well-defined smooth morphism $\psi \colon \tilde{\bbP }^N_{B_1}\to \bbP ^{n-1}_{B_1}$.
\begin{equation*}
\xymatrix@R=7pt{
&\tilde{X}_1 \ar[dd]^{\pi }	\ar@{^{(}->}[r]\ar@{}|{|}[r]
&\tilde{\bbP }^N_{B_1} \ar[dd] \ar[r]_\psi
&\bbP ^{n-1}_{B_1}
\\%
X_2\ar@{^{(}->}[dr]\ar@{}|{\circ }[dr]
\ar@{^{(}->}[ur]\ar@{}|{\circ }[ur]
&&&
\\%
&\bar{X}_1 		\ar@{^{(}->}[r]\ar@{}|{|}[r]
&{\bbP }^N_{B_1}
&
}
 \end{equation*}

Some observations on this map are in order.
Given a point 
%%
%\begin{equation*}
$y=[\eta _1:\dots :\eta _n ]\colon $ $\Spec (\Omega )\to \bbP ^{n-1}_{B_1}$
% \end{equation*}
%%
($\Omega $ is a field),
the fiber of $\psi $ at $y$ is identified with the linear subscheme $V(y)$ of $\bbP ^N_{\Omega }$ defined by the equations
$\eta _j \psi _k(X_0,\dots ,X_N)=\eta _k \psi _j (X_0,\dots ,X_N)$
($j,k\in \{ 1,\dots ,n \} $).
If we choose an index $j$ with $\eta _j\neq 0$, a minimal set of equations can be chosen as 
\begin{equation}\label{Eq:linear-y}
\eta _j \psi _k(X_0,\dots ,X_N)=\eta _k \psi _j (X_0,\dots ,X_N)
\qquad k\in \{ 1,\dots ,n \} \setminus \{ j \} .
 \end{equation}
It contains $V(\psi )_\Omega $ and isomorphic to $\bbP ^{N-n+1}_\Omega $.

Let $E(\psi )\subset \tilde{\bbP }^{N}_{B_1}$ be the exceptional divisor and $\pi \colon E(\psi )\to V(\psi )$ the restriction of the blow-up map.
\begin{equation*}
\xymatrix@R=6pt@C=6pt{
&E(\psi )\ar[dl]_{\pi }\ar[dr]^{\psi }&
\\%
V(\psi )\ar[dr]&&\bbP ^{n-1}_{B_1}\ar[dl]
\\%
&B_1 &
}
 \end{equation*}
It turns out that this square is cartesian.

Now, let us observe that the restriction $\psi _{|\tilde{X}_{1}}\colon \tilde{X}_{1}\to \bbP ^{n-1}_{B_1}$ is equi-dimensional of fiber dimension $1$.
Indeed, given a point $y\in \bbP ^{n-1}_{B_1}(\Omega ) $, consider $(\bar{X}_1)_\Omega := \bar{X}_1\times _{B_1} \Omega $ $\subset \bbP ^N_\Omega $ which has dimension $n $.
The fiber of $\psi _{|\tilde{X}_1}$ over $y$ is contained in the intersection $(\bar{X}_1)_\Omega \cap V(y)$ in $\bbP ^N_{\Omega }$.
We know that $V(y)$ contains $V(\psi )_\Omega $ as a linear subscheme of codimension $1$.
Every irreducible component $Y_\alpha $ of $ (\bar{X}_1)_\Omega \cap V(y)$ meets $V(\psi )_\Omega $ because it is the intersection in a projective space, and we chose $\psi $ so that the intersection $Y_\alpha \cap V(\psi )_{\Omega }$ has dimension $0$. %By Krull's Hauptidealsatz,
Therefore we conclude that $Y_\alpha $ has dimension $1$.

\ \ \

Set $T:= \bbP ^{n-1}_{B_1}$ and $t:= \psi (x)\in \bbP ^{n-1}_{B_1}$.
Choose any projective embedding $\tilde{X}_1\hookrightarrow \bbP ^{N_2}_T$.
Set $(\tilde{X}_1)_t:= \tilde{X}_1\times _T t$ which is in $\bbP ^{N_2}_t$.
Also, set $(X_2)_t:=X_2\times _T t$ which is an open subset of $(\tilde{X}_1)_t$.
There is a hypersurface $H_t$ in $\bbP^{N_2}_t$ satisfying the next three conditions. [Indeed, there exist hypersurfaces satisfying (i) in some high degree; a general one among such (of any fixed degree) satisfies (ii) and (iii).]

\begin{enumerate}[(i)]
\item {}[In case $x$ is a closed point of $(X_2)_t$:]
$H_t$ contains $x$.

\item $(\tilde{X}_1)_t$ and $H_t$ meet properly in $\bbP^{N_2}_t$.

\item {}[Let $(X_2)^-_t$ be the closure of $(X_2)_t$ in $\bbP^{N_2}_t$. Then $(X_2)^-_t\setminus (X_2)_t$ is a finite set. That being said:]
$H_t$ does not meet $(X_2)^-_t\setminus (X_2)_t$.

 \end{enumerate}

Let $T^{\hen }\to T$ be the henselization of $T$ at $t$.
Let $H_{T^\hen }\subset \bbP ^{N_2}_{T^\hen }$ be any flat family of hypersurfaces that specializes to $H_t$.
Base change by $T^{\hen }\to T$ will be denoted by the subscript $(-)_{T^\hen }$.
The closed subscheme $D_{T^\hen }:=H_{T^\hen }\times _{\bbP ^{N_2}_{T}} \tilde{X}_1$ of $(\tilde{X}_1)_{T^\hen }$ is an effective Cartier divisor
because $(\tilde{X}_1)_{T^\hen }$ is reduced (hence without embedded components).
By condition (ii), the scheme $D_{T^\hen }$ is the disjoint union of its local components finite over $T^\hen $; in particular, it decomposes as
\begin{equation*}
D_{T^\hen }=\mcal{D}_{T^\hen }\sqcup \mcal{D}'_{T^\hen }
 \end{equation*}
where $\mcal{D}_{T^\hen }$ has closed points in the open set $(X_2)_{t }$ (and hence it is contained in $(X_2)_{T^\hen }$) and $\mcal{D}'_{T^\hen }$ is disjoint from $(X_2)_t$.

By the usual limit argument as in \cite[IV${}_2$ \S 8]{EGA},
there is an affine Nisnevich neighborhood $T_3\to T$ of $t$ such that: we can take a hypersurface $H_{T_3}$ such that $H_{T^\hen }=H_{T_3}\times _{T_3}T^\hen $; the decomposition of $D_{T^\hen }$ comes from a decomposition $D_{T_3}=\mcal{D}\sqcup \mcal{D}'$; and the conditions stated in the previous paragraph become true already after the base change $T_3\to T$. Let us denote $\tilde{X}_3:=\tilde{X}_1\times _T T_3$ and
$X_3:=X_2\times _T T_3$. Then we have:
\begin{equation*}
\xymatrix{
X_3 \ar@{^{(}->}[r]\ar|{\circ}[r] &\tilde{X}_3\ar[rr]^{\text{projective,}}_{\text{1-dim}}&&T_3\ar[r]^{\text{Nis}}&\bbP ^{n-1}_{B_1}
\\%
&\mcal{D} \ar@{^{(}->}[ul]\ar|{\diagup}[ul]^{\text{Cartier div.}} \ar@{^{(}->}[u]\ar|{-}[u] \ar[rru]_{\text{finite}}
}
 \end{equation*}
For integers $m$, we can consider the line bundles $\mcal{O}_{\tilde{X}_3}(m\mcal{D})$ on $\tilde{X}_3$. We claim that for any sufficiently large $m>0$, the restriction map
\begin{equation*}
\Gamma (\tilde{X}_3,\mcal{O}_{\tilde{X}_3}(m\mcal{D}))
\to \Gamma (\mcal{D},\mcal{O}_{\mcal{D}}(m\mcal{D}))
 \end{equation*}
is surjective.
Indeed, we have an exact sequence
$0\to \mcal{O}_{\tilde{X}_3}(-\mcal{D})\to \mcal{O}_{\tilde{X}_3}\to \mcal{O}_{\mcal{D}}\to 0$ by definitions, which gives exact sequences
\begin{equation}\label{Eq:short-exact}
0\to \mcal{O}_{\tilde{X}_3}(m\mcal{D})\to \mcal{O}_{\tilde{X}_3}((m+1)\mcal{D})\to \mcal{O}_{\mcal{D}}((m+1)\mcal{D})\to 0.
 \end{equation}
Since $\mcal{D}$ is an affine scheme, we have
$H^1(\mcal{D},\mcal{O}_{\mcal{D}}((m+1)\mcal{D}))=0$.
Therefore we have surjections
\begin{align*}
H^1(\tilde{X}_3,\mcal{O}_{\tilde{X}_3}(m\mcal{D}) )
&\twoheadrightarrow
H^1(\tilde{X}_3,\mcal{O}_{\tilde{X}_3}((m+1)\mcal{D}) )
\twoheadrightarrow \dots 
\\%
&\twoheadrightarrow
H^1(\tilde{X}_3,\mcal{O}_{\tilde{X}_3}((m+k)\mcal{D}) )
\twoheadrightarrow \dots .
 \end{align*}
Since $M_m:= H^1(\tilde{X}_3,\mcal{O}_{\tilde{X}_3}(m\mcal{D}) )$ is a finitely generated $\Gamma (T_3,\mcal{O}_{T_3})$-module, which is noetherian, the group $\ker (M_m\to M_{m+k})$ stabilizes as $k$ tends to $\infty $. Therefore $M_{m+k}\to M_{m+k+1}$ is an isomorphism for any sufficiently large $k$.
It follows that in the cohomology long exact sequence associated to the short exact sequence \eqref{Eq:short-exact} with the index shifted:
\begin{equation*}
H^0(\tilde{X}_3,\mcal{O}_{\tilde{X}_3}(m\mcal{D}))\to
H^0(\mcal{D},\mcal{O}_{\mcal{D}}(m\mcal{D}))
\to M_{m-1} \to M_{m},
 \end{equation*}
the last map is bijective for any sufficiently large $m$.
We conclude that the map
$H^0(\tilde{X}_3,\mcal{O}_{\tilde{X}_3}(m\mcal{D}))\to
H^0(\mcal{D},\mcal{O}_{\mcal{D}}(m\mcal{D}))$
is surjective.

By this surjectivity, we can find a section
$s_0\in \Gamma (\tilde{X}_3,\mcal{O}_{\tilde{X}_3}(m\mcal{D}))$
which maps to a nowhere vanishing section of 
$\Gamma (\mcal{D},\mcal{O}_{\mcal{D}}(m\mcal{D}))\cong \Gamma (\mcal{D},\mcal{O}_{\mcal{D}})$.
Let $s_1\colon \mcal{O}_{\tilde{X}_3}\to \mcal{O}_{\tilde{X}_3}(m\mcal{D})$ be the canonical inclusion. Since the zero-loci of $s_0$ and $s_1$ are disjoint (one is away from $\mcal{D}$ and the other is $|\mcal{D}|$),
we have a well-defined morphism
\begin{equation*}
f:= (s_0:s_1)\colon \tilde{X}_3\to \bbP^1_{T_3}.
 \end{equation*}
Let us consider its restriction to $(\tilde{X}_3)_t=(\tilde{X}_1)_t$
which is a projective curve over $k(t)$.
The inverse image of $\infty \in \bbP ^1_t$ is exactly 
$(\tilde{X}_1)_t\cap \mcal{D}$.
Let $Y_\alpha $ be any irreducible component of $(\tilde{X}_1)_t$.
If it meets the open subset $(X_3)_t=(X_2)_t$, then $Y_\alpha \cap \mcal{D}$ is a non-empty finite set by the definition of $\mcal{D}$ (it is non-empty by condition (iii)); in this case $Y_\alpha \to \bbP ^1_t$ is finite flat.
To the contrary, if $Y_\alpha $ is away from $(X_3)_t$, we have $Y_\alpha \cap \mcal{D}=\emptyset $; in this case $Y_\alpha \to \bbP ^1_t$ is a constant map to a point.

The quasi-finite locus of a morphism is open \cite[IV${}_3$ 13.1.4]{EGA}.
We can find an affine open neighborhood $T_4\subset T_3$ of $t$ such that
$\mcal{D}_4:=\mcal{D}_{T_4}$ is contained in the quasi-finite locus of $f$.
Let $f_4\colon \tilde{X}_4\to \bbP ^1_{T_4}$ be the base change $f\times _{T_3}T_4$ and let $X_4\subset X_3$ be the quasi-finite locus of $f_{|X_3\times _{T_3} T_4}$.
\begin{equation*}
\xymatrix{
f^{-1}(\infty _{T_4})=\mcal{D}_4\ar@<25pt>[d]\ar@{^{(}->}[r]\ar@{}|(0.7){|}[r]
&X_4 \ar[d]_{\begin{subarray}{r}\text{quasi-}\\ \text{finite}\end{subarray} } 
\ar@{^{(}->}[r]\ar@{}|{\circ}[r]
& \tilde{X}_4\ar[dl]^f
\\%
{}\phantom{f^{-1}(\infty _{T_4})=}\infty _{T_4}\ar@{^{(}->}[r]\ar@{}|(0.7){|}[r] 
&\bbP^1_{T_4}
}
 \end{equation*}
Then the subset $W_4:= f(\tilde{X}_4\setminus X_4)\subset \bbP ^1_{T_4}$
is proper over $T_4$ and is contained in $\bbP ^1_{T_4}\setminus \infty _{T_4} =\bbA ^1_{T_4}$; therefore it is finite over $T_4$.
The induced map from $\tilde{X}_4\setminus f^{-1}(W_4)$ to $\bbP ^1_{T_4}\setminus W_4$ is finite because it is proper and quasi-finite (as $f$ is quasi-finite on $X_4$).
%Set $X_5:= \tilde{X}_4\setminus f^{-1}(W)$. It 
Note that $\tilde{X}_4\setminus f^{-1}(W_4)$ contains $x$ by condition (i).

\ \ \

Now, by induction on $n$, there are Nisnevich neighborhoods $B_5\to B_1$ of $b$ and $T_5\to T_4\times _{B_1}B_5$ of $t$ such that there is a projective compactification $T_5\hookrightarrow \bar{T}_5$ which is fiberwise dense over $B_5$.
Base change of $T_4$-schemes by $T_5\to T_4$ will be indicated by replacement of the number.
%Replacement of ``4'' by ``5'' in our construction so far shall mean the base change $T_5\to T_4$.

Take any factorization of $f$ of the form $\tilde{X}_5\climm \bbP ^{N_5}_{T_5} \times _{T_5}\bbP ^1_{T_5}  \to \bbP^1_{T_5}$ and denote by 
$\bar{X}_5$ the reduced closure of $\tilde{X}_5$ in $\bbP ^{N_5}_{\bar{T}_5} \times _{\bar{T}_5}\bbP ^1_{\bar{T}_5}$.
We have the following diagram where every square is cartesian:
\begin{equation*}
\xymatrix{
X_6:=\tilde{X}_5\setminus f\inv (W_5)
\ar@{^{(}->}[r]\ar@{}|(0.7){\circ}[r]
\ar@{^{(}->}[d]\ar@{}|{-}[d]
&\tilde{X}_5
\ar@{^{(}->}[r]\ar@{}|{\circ}[r]
\ar@{^{(}->}[d]\ar@{}|{-}[d]
&\bar{X}_5
\ar@{^{(}->}[d]\ar@{}|{-}[d]
\\%
\bbP ^{N_5}_{{T}_5} \times _{{T}_5}(\bbP ^1_{{T}_5}\setminus W_5)
\ar@{^{(}->}[r]\ar@{}|(0.55){\circ}[r]
\ar[d]
&\bbP ^{N_5}_{{T}_5} \times _{{T}_5}\bbP ^1_{{T}_5}
\ar@{^{(}->}[r]\ar@{}|{\circ}[r]
\ar[d]
&\bbP ^{N_5}_{\bar{T}_5} \times _{\bar{T}_5}\bbP ^1_{\bar{T}_5}
%\ar[r]
\ar[d]
%&\bbP ^1_{\bar{T}_5}
\\%
\bbP ^1_{{T}_5}\setminus W_5
\ar@{^{(}->}[r]\ar@{}|{\circ}[r]
&\bbP ^1_{T_5}
\ar@{^{(}->}[r]\ar@{}|{\circ}[r]
&\bbP ^1_{\bar{T}_5}
} \end{equation*}
Take the Stein factorization \cite[III${}_1$ \S 4.3]{EGA} of the composite map $\bar{f}_5\colon \bar{X}_5\to \bbP ^1_{\bar{T}_5}$:
\begin{equation*}
\bar{X}_5\to \bar{X}_6\xrightarrow{\text{finite}} \bbP ^1_{\bar{T}_5}
 \end{equation*}
so that the first map $\bar{X}_5\to \bar{X}_6$ has geometrically connected fibers. Since $\bar{f}_5$ is already finite over the open set $\bbP^1_{T_5}\setminus W_5$ of the target, we know $\bar{X}_6\times _{\bbP^1_{\bar{T}_5}} (\bbP ^1_{T_5}\setminus W_5)$
is isomorphic to $X_6:=\tilde{X}_5\setminus f\inv (W_5)$.

\ \ \

Now we claim that the open immersion $X_6\opimm \bar{X}_6$ is fiberwise dense over $B_5$.
First, note that $\bar{T}_5$ and $\bar{X}_6$ are equi-dimensional over $B_5$ with fiber dimensions $n-1$ and $n$ by the same reasoning as in the second paragraph of this proof (using the fact that $B_5$ is Dedekind).
Let $b'\in B_5$ be any point and $Y_\beta \subset (\bar{X}_6)_{b'}$ be any irreducible component of the fiber. We have to prove that $Y_\beta $ meets the open set $X_6$.

The composition $Y_\beta \climm (\bar{X}_6)_{b'}\to (\bbP ^1_{\bar{T}_5})_{b'}$
is a finite morphism of $n$-dimensional $k(b')$-schemes, and so it is a surjection to an irreducible component.
The open set $(\bbP ^1_{T_5}\setminus W_5)_{b'}$ of the target is dense because
$(T_5)_{b'}\subset (\bar{T}_5)_{b'}$ is dense and $W_5$ is finite over $T_5$.
Therefore its inverse image to $Y_\beta $, which is exactly $Y_\beta \cap X_6$, is non-empty.
This verifies that $X_6\subset \bar{X}_6$ is fiberwise dense over $B_5$.

\ \ \

It remains to find claimed affine and projective embeddings.
The scheme $\bar{X}_6$ is projective over $B_5$ because the morphisms $\bar{X}_6\to \bar{T}_5$ and $\bar{T}_5\to \bar{B}_5$ are, so take any closed immersion $\bar{X}_6\climm \bbP ^{N_6}_{B_5}$.
Since the ideal of the closed subscheme $(\bar{X}_6\setminus X_6)_{\mrm{red}}$
is generated by global sections after sufficient Serre twist, there exists a flat family of hypersurfaces $H_{B_5}\subset \bbP ^{N_6}_{B_5}$
containing $\bar{X}_6\setminus X_6$ and such that its fiber $H_b$ over $b$ meets $(\bar{X}_6)_b$ properly in $\bbP^{N_5}_{b}$.
Also, we require that $H_b$ does not contain $x$.
There is an affine open neighborhood $B_7\to B_5$ of $b$ such that the base change $H_{B_7}$ mees $\bar{X}_{7}:=\bar{X}_6\times _{B_5}B_7$ fiberwise properly over $B_7$.

Let $d$ be the degree of $H_{B_7}$.
Then via the $d$-fold Veronese embedding $\bbP ^{N_5}_{B_7}\climm \bbP ^{N_7}_{B_7}$ followed by a linear automorphism, $H_{B_7}$ is the restriction of the infinity-hyperplane of $\bbP^{N_7}_{B_7}$.
Put $X_7:= \bar{X}_7\setminus H_{B_7}$ and consider $\bar{X}_7$ as a closed subscheme of $\bbP ^{N_7}_{B_7} $.
Then the cartesian diagram
\begin{equation*}
\xymatrix{
X_7 \ar@{^{(}->}[r]\ar@{}|{\circ }[r]
\ar@{^{(}->}[d]\ar@{}|{-}[d]
& \bar{X}_7 \ar@{^{(}->}[d]\ar@{}|{-}[d]
\\%
\bbA ^{N_7}_{B_7} \ar@{^{(}->}[r]\ar@{}|{\circ }[r] & \bbP ^{N_7}_{B_7}
}
\end{equation*}
verifies our assertion.
\end{proof}
%%% Proof ends %%%

%%%%%%%%%%%%%%%%%%%%%%%%%%%%%%%
%%%%%%%%%%%%%%%%%%%%%%%%%%%%%%%
\section{General case}\label{Sec:General-case}

%%%%%%%%%%%%%
%% Theorem %%
\begin{definition}
Let $\mcal{C}$ be a family of constructible subsets of an equi-dimensional scheme $X$ and $q\ge 0$ be an integer. Define the {\em shifted} family $\mcal{C}[q]$ by
$(\mcal{C}[q])_d:=\mcal{C}_{d-q}$ if $d<\dim (X)$ and $:= X$ if $d\ge \dim (X)$.

For example $\mcal{C}[1]$ is the following (write $n:=\dim (X)$): 
\begin{equation*}\xymatrix@R=0pt@C=10pt{
&-1&0&1&2&\dots &\dim (X)-1 &\dim (X)
\\%
\mcal{C}\phan{[1]}: &\emptyset \ar@{}|{\subseteq }[r] & C_0 \ar@{}|{\subseteq }[r]&C_1 \ar@{}|{\subseteq }[r]&C_2\ar@{}|{\subseteq }[r]&\dots \ar@{}|{\subseteq }[r]&C_{n-1} \ar@{}|{\subseteq }[r]&X
\\%
\mcal{C}[1]: &\emptyset \ar@{}|{\subseteq }[r]&\emptyset \ar@{}|{\subseteq }[r]& C_0 \ar@{}|{\subseteq }[r] &C_1 \ar@{}|{\subseteq}[r] &\dots \ar@{}|{\subseteq }[r]&C_{n-2} \ar@{}|{\subseteq }[r]&X
}
\end{equation*}
Note that we have a chain law $(\mcal{C}[p])[q]=\mcal{C}[p+q]$ and $\mcal{C}[n]$ is the trivial family.
\end{definition}
%% Theorem ends %%
%%%%%%%%%%%%%%%%%%

We go back to the notation of Main Theorem \ref{Th:full-statement}; namely, we are given an effective Cartier divisor $D\subset X$ of an equi-dimensional $k$-scheme of dimension $n$ and a family $\mcal{C}$ of constructible subsets of $X\setminus D$, and an equi-dimensional $k$-scheme $Y$.
To prove that the inclusion $z^{i}_{\mcal{C}}(-\times Y\mid D_-\times Y,\bullet )\subset z^{i}(-\times Y\mid D_-\times Y,\bullet )$ is a quasi-isomorphism on $X_\Nis $, it suffices to show that the quotient complex $z^{i}_{\mcal{C}[q]}(-\times Y\mid D_-\times Y,\bullet )/z^i_{\mcal{C}[q-1]}
$
is Nisnevich locally acyclic for $1\le q \le n$, for then in the sequence below:
\begin{equation*}
z^{i}_{\mcal{C}}\subset z^{i}_{\mcal{C}[1]}\subset \dots
z^{i}_{\mcal{C}[n]}=z^{i}(-\times Y\mid D_-\times Y,\bullet ),
 \end{equation*}
every inclusion is a Nisnevich local quasi-isomorphism on $X$.
Refreshing the notation $\mcal{C}:= \mcal{C}[q-1]$, it suffices to prove that
the quotient complex $z^{i}_{\mcal{C}[1]}(-\times Y\mid D_-\times Y,\bullet )/z^{i}_{\mcal{C}}$ is locally acyclic.

Let us recall the following principle: To show that a given complex $Z_\bullet $ of abelian groups is acyclic, it suffices to show that for every finitely generated subcomplex $Z_\bullet '$ (i.e., it has only finitely many non-trivial terms, each of which is a finitely generated abelian group), the inclusion $Z_\bullet '\subset Z_\bullet $ induces the zero map on homology.
This is because every homology class of $Z_\bullet $ comes from some $Z_\bullet '$.

\begin{proof}[Proof of Theorem \ref{Th:full-statement}]

Since the assertion is local on $X$, we may assume that $X$ is affine and $D$ is a principal divisor defined by a function $\pi \colon X\to \bbA ^1$.
We take an arbitrary finite set 
$\{ V_\lambda \subset $
$X\times Y\times \square ^{p(\lambda )}\} _\lambda $
of irreducible cycles
belonging to $z^{i}_{\mcal{C}[1]}(X\times Y\mid D\times Y,\bullet )$ which is closed under the passage from a cycle $V_\lambda $ to the irreducible components of its restriction to faces.

Consider the subcomplex
\begin{equation*}
\frac{z^{i\prime }_{\mcal{C}[1]}}{z^{i}_{\mcal{C}}}
\quad\text{ of }\quad \frac{z^{i}_{\mcal{C}[1]}(X\times Y|D\times Y,\bullet )}{z^{i}_{\mcal{C}}(X\times Y|D\times Y,\bullet ) }
\end{equation*}
generated by the cycles $\{ V_\lambda \} _\lambda $;
we have to find a Nisnevich neighborhood $X'\to X$ of a given point $x\in X$ such that the map (where $D':=D\times _X X'$)
\begin{equation*}
\frac{z^{i\prime }_{\mcal{C}[1]}}{z^{i}_{\mcal{C} }} \to \frac{z^{i}_{\mcal{C}[1]}(X'\times Y|D'\times Y,\bullet )}{z^{i}_{\mcal{C} }(X'\times Y|D'\times Y,\bullet )}
 \end{equation*}
induces the zero map on homology.

Consider all the subsets of the form
$V_\lambda \cap (C_d \times Y\times \square ^{p(\lambda )})$,
and write $S\subset X\setminus D$ for the finite set consisting of the projections of their generic points to $X$.

%%%%%%%%%%%%%
%% Theorem %%
\begin{lemma}\label{Lem:Smooth-replacement}
There exist an affine open neighborhood $X_1$ of $x$ in $X$ and a unit $f\in \mcal{O}(X_1)^*$ such that the function $f\pi $: $X_1\to \bbA ^1$ is smooth on the finite set $S\cap X_1$.
\end{lemma}
%% Theorem ends %%
%%%%%%%%%%%%%%%%%%

[The following proof would become simpler if $x$ is in the smooth locus of $X$.]

Let $X\climm \bbA ^N_k $ be any closed immersion and consider linear projections $\phi \colon \bbA ^N\to \bbA ^1$. For a general $\phi $,
the induced map $f_0:= \phi _{|X} \colon X\to \bbA ^1_k$ is smooth at the points of $S$ (e.g.\ by \S \ref{Sec:smoothness}).
Fix such an $f_0$. We are going to choose $f$ in the form 
\begin{equation*}
f=f_0+a \qquad
(a\in k). 
 \end{equation*}
The morphism $f\pi $ is smooth at a given point $x_i\in X$
if and only if the induced $k(x_i)$-linear map
[let $t$ be the coordinate for $\bbA ^1_k$]
\begin{equation*}
(f\pi )^*\colon k(x_i)\cdot dt \to T_{x_i}X
 \end{equation*}
is injective. Here, the basis $dt$ is mapped to the next element of
$T_{x_i}X$:
\begin{align*}
d(f\pi )(x_i) &= \pi (x_i) \cdot (df)(x_i) +f(x_i)\cdot (d\pi )(x_i)
\\%
	& = \pi (x_i)\cdot (df_0)(x_i)+(f_0(x_i)+a) \cdot (d\pi )(x_i)
	 .
 \end{align*}
Since the first term on the right hand side is non-zero if $x_i\in S$, this element is non-zero for all but possibly one value of $a\in k$.
Since $S$ is a finite set, these elements are simultaneously non-zero for all $x_i\in S$ if $a\in k$ avoids a certain finite set.

Lastly, the function $f_0+a $ is invertible at each point of $S\cup \{ x\} $ for a general $a\in k$. 
For such a general choice of $a\in k$, set $X_1:= f^{-1}(\bbA ^1\setminus \{ 0\} )$.
This completes the proof of Lemma \ref{Lem:Smooth-replacement}.

\ \ \

Returning to the proof of Theorem \ref{Th:full-statement}, let $B\to \bbA ^1_k$ be the henselization at the origin.
We apply Theorem \ref{Th:Noether} to the equi-dimensional morphism
$f\pi \colon X_1\to \bbA ^1_k$ of fiber dimension $n-1$
to find a Nisnevich neighborhood $X'\to X_1\times _{\bbA ^1} B$ of $x$ 
and an embedding $X'\climm \bbA ^{N}_B$ such that the following is true: If $\bar{X}'\subset \bbP ^{N}_B$ is the closure of $X'$ in the projective space, then its intersection with the infinity-hyperplane $\bar{X}'\cap (\bbP ^{N}_{B }\setminus \bbA ^N_B)$ is equi-dimensional over $B$ of fiber dimension $n-2$.
Base-change from $X$ to $X'$ shall be indicated by primes $(-)'$: For example $D':=D\times _X X'$ and $C_d'=C_d\times _X X'$. 

It follows that there is a dense open subset $U$ of $(M_{N\times (n-1)})_{\{ 0\} }$ such that 
if a matrix $\phi \in M_{N\times (n-1)}(B)$ specializes %at the origin 
to a matrix in $U$, then the induced map $\phi _{|X'}\colon X'\to \bbA ^{n-1}_B$ is finite and surjective (hence flat on the smooth locus of $X'$ \cite[IV${}_{2}$ 6.1.5]{EGA}).
In particular, pull-back and push-forward by $\phi _{|X'}$ preserve the modulus condition of cycles, where $X' $ is equipped with the modulus $D'$ and $\bbA ^{n-1}_B$ with $\bbA ^{n-1}_{\{ 0\} }$.
For such a $\phi $ and each of $V_\lambda $,
%\in z^{i }(X\times Y|D\times Y,p(\lambda ))$, 
we may consider the cycle 
$(\phi _{|X'}^*\phi _{|X'\ *}(V_\lambda ')) -(V_\lambda ')$
which is an element of $z^i(X'\times Y|D'\times Y,p(\lambda ))$.
We denote it by $\phi ^*\phi _*V_\lambda '-V_\lambda '$ for short. 
Below, when we say that a general $\phi \in M^*_{N\times (n-1)}$ has some property, we mean that there is a dense open subset $U'$ of 
$(M^*_{N\times (n-1)})_{\Frac (B)}$ such that every $\phi \in U'$ specializing into $U$ has the property.
(In this case, such a $\phi $ {\it exists} because the matrices in $(M^*_{N\times (n-1)})_{\Frac (B)}$ which specialize into $U$ form a dense subset.) 

\ \ \

By the usual limit argument, the cycle complex $z^i(X'\times Y|D'\times Y,\bullet ) $ is the direct limit of complexes $z^i(X_\nu \times Y|D_\nu \times Y,\bullet )$ where $X_\nu \to X$ runs through certain Nisnevich neighborhoods of $x$ of finite type.
For a family $\mcal{C}$ of constructible subsets of $X\setminus D$, we write $z^i_{\mcal{C}}(X'\times Y|D'\times Y,\bullet )$ for the limit of the subcomplexes $z^i_{\mcal{C}}(X_\nu \times Y|D_\nu \times Y,\bullet )$ (we are omitting the pull-back notation for $\mcal{C}$ as usual).
In fact, the literal definition of that group does not make sense because $X'\setminus D'$ has dimension $n-1$.
%, but we choose to write this way.
The interested reader will find that 
%if we define $\mcal{C}'=\{ C_d' \} _{d\in \bbZ } $
%by $C_d':=C_{d+1}\times _X X'$,
%for $d< n-1$ and $C_d'= X'\setminus D'$ for $d\ge n-1$,
%then 
our $z^i_{\mcal{C}}(X'\times Y\mid D'\times Y,\bullet )$ is equal to the literally defined $z^i_{\mcal{C}[-1]'}(X'\times Y\mid D'\times Y,\bullet )$.
We adopt the same convention with $z^i_{\phi _*\mcal{C}}(\bbA ^{n-1}_B\times Y|\bbA ^{n-1}_{\{ 0\} }\times Y ,\bullet )$ appearing below.
%When we make statements about dimensions of subsets of $X'\setminus D'$,
%they are to be understood as statements about corresponding subsets in some $X_\nu \setminus D_\nu $ for a large enough index $\nu $.

\ \ \

By \S \ref{Sec:birationality} applied to the scheme $X'\setminus D'$ over the field $\Frac (B)$, we know that a general $\phi \in M^*_{N\times (n-1)}$ has the property that
the support of the cycle $\phi ^*\phi _*V_\lambda '-V_\lambda '$ is $\phi ^+(V_\lambda ')$.

%Since we have modified the defining equation for $D$ using 
Thanks to
Lemma \ref{Lem:Smooth-replacement}, we may apply \S \ref{Sec:Chow's-lemma} to subsets $V_\lambda '$ and $W:= C_d'$ of $X'\setminus D'$.
It follows that a general $\phi \in M^*_{N\times (n-1)}$ satisfies the following 
for all $\lambda $ and $d$:
\begin{align*}
\dim &(|\phi ^*\phi _* V_\lambda '-V_\lambda '|\cap (C_d'\times Y\times \square ^{p(\lambda )}))
\\%
&\le \max \{  \dim (V_\lambda ' \cap (C_d'\times Y\times \square ^{p(\lambda )}))-1,
\dim (V_\lambda ')+\dim (C_d')-\dim (X')  \} .
 \end{align*}
This formula precisely means that $\phi ^*\phi _*V_\lambda '-V_\lambda '$ is in the smaller subcomplex 
$z^i_{\mcal{C} }(X'\times Y|D'\times Y,\bullet )$.
In particular, we deduce that $\phi ^*\phi _*V_\lambda '$ is in $z^i_{\mcal{C}[1]}(X'\times Y|D'\times Y,\bullet )$.
Hence, we have the equality
\begin{equation*}
\phi ^*\phi _*V_\lambda ' = V_\lambda '
\quad \text{ in }\quad
\frac{z^{i }_{\mcal{C}[1]}(X'\times Y|D'\times Y,\bullet )}{z^i_{\mcal{C} }(X'\times Y|D'\times Y,\bullet )}.
 \end{equation*}
It also follows that
$\phi _*V_\lambda '$ is in $z^i_{\phi _*\mcal{C}[1]}(\bbA ^{n-1}_B\times Y|\bbA ^{n-1}_{\{ 0\} }\times Y,\bullet )$.
This is because the subset $(\phi _*V_\lambda ')\cap (\phi (C_d) \times Y\times \square ^{p(\lambda )})$ is equal to $\phi (\phi ^*\phi _*V_\lambda ' \cap (C_d\times Y\times \square ^{p(\lambda )}))$ by the usual `projection formula' of subsets.
If $V_\lambda $ already happens to be in the smaller subcomplex $z^i_{\mcal{C} }(X\times Y|D\times Y,\bullet )$,
the same reasoning gives: 
\begin{align*}
&\phi ^*\phi _*V_\lambda ' \in z^i_{\mcal{C} }(X'\times Y|D'\times Y,\bullet ) \\%
\text{and }\quad &\phi _*V_\lambda '\in z^i_{\phi _*\mcal{C} }(\bbA ^{n-1}_B\times Y|\bbA ^{n-1}_{\{ 0\} }\times Y,\bullet ).
 \end{align*}
Therefore, the first map of the next sequence is well-defined:
\begin{equation*}
\frac{z^{i\prime }_{\mcal{C}[1]}}{z^i_{\mcal{C} }}
\xrightarrow{\phi _*} \frac{z^{i }_{\phi _*\mcal{C}[1]}(\bbA ^{n-1}_B\times Y|\bbA ^{n-1}_{\{ 0\} }\times Y,\bullet )}{z^i_{\phi _*\mcal{C} }}
\xrightarrow{\phi ^*}
\frac{z^{i }_{\mcal{C}[1]}(X'\times Y|D'\times Y,\bullet )}{z^i_{\mcal{C} }}
 \end{equation*}
and the composite is equal to the pull-back map by $X'\to X$.
On the other hand, since the middle complex is acyclic by Theorem \ref{Th:affine-case},
the composite induces the zero map on homology.
This completes the proof of Theorem \ref{Th:full-statement}.
\end{proof}

For Theorem \ref{Th:additive-Chow} (the case of additive higher Chow),
we proceed similarly using affine embedding $X\climm \bbA ^N_k$ and general linear projection
$\bbA ^N_k\surj \bbA ^{\dim (X)}_k $,
and use the fact that the assertion is known for the pair $(\bbA ^{\dim (X)+n},\bbA ^{\dim (X)}\times D)$
by Theorem \ref{Th:affine-additive-case}.
For more details, the reader can consult 
\cite[\S 6]{KrishnaParkAdditive}.

%\input{TexSimplicial}

%%%%%%%%%%%%%%%%%%%%%%%%%%%%%%%
%%%%%%%%%%%%%%%%%%%%%%%%%%%%%%%
\section{Remarks on the simplicial version}\label{Sec:simplicial}

Here we explain briefly that the content of this paper is true for the (obvious) simplicial analog of the Binda-Saito cycle complex.

Also, in Remark \ref{Rem:cheaper-simplicial}, we briefly mention Miyazaki's result \cite{MiyazakiPrivate} that the simplicial version is equivalent to the cubical, at least in the pro setting along the multiples of the divisor $D\subset X$.

The simplicial version could be important because, 
in addition to the fact that the simplicial terminology is much more preferred in general,
it turns out that the Chern classes with modulus from the relative $K$-groups $K_{\ge 0}(X,D)$ to the Nisnevich hypercohomology of the cycle complex factors through the simplicial version \cite[footnote in \S 2.2]{IwasaKai}.

%%%%%%%%%%%%%%%%%%%%%%%%%%%%%%%
%%%%%%%%%%%%%%%%%%%%%%%%%%%%%%%
\subsection{Definition}

Write $\Delta ^p:= \Spec (\bbZ [t_0,\dots , t_p]/(\sum _{i=0}^p t_i -1))$
and $\bar{\Delta }^p:= \Proj (\bbZ [U,T_0,\dots ,T_p]/(-U+\sum _{i=0}^pT_i))$
so that $\Delta ^p$ is the non-vanishing locus of $U$ in $\bar{\Delta }^p$.
Let $\bar{\Delta }^p_\infty $ be the effective Cartier divisor in $\bar{\Delta }^p$ defined by $U=0$.
When we work over a base ring, we consider them over the ring.

We now define a simplicial variant of the modulus condition.
Let $(X,D)$ be a pair of an algebraic scheme and an effective Cartier divisor.

%%%%%%%%%%%%%
%% Theorem %%
\begin{definition}
A closed subset $V$ of $X\times \Delta ^p $ is said to satisfy the modulus condition if the following is true:
Let $\bar{V}^N$ be the normalization of the (reduced) closure $\bar{V} $ of $V$ in $X\times \bar{\Delta }^p$.
Then the inequality of Cartier divisors on $\bar{V}^N$
\begin{equation*}
D|_{\bar{V}^n}\le \bar{\Delta }^p_\infty |_{\bar{V}^n}
 \end{equation*}
holds, where $(-)|_{\bar{V}^n}$ denotes the pull-back to $\bar{V}^N$.

\end{definition}
%% Theorem ends %%
%%%%%%%%%%%%%%%%%%

Pull-back along the structure morphisms $\Delta ^{q}\to \Delta ^p$
of the cosimplicial scheme $\{ \Delta ^p\} _p $ turns out to preserve the modulus condition:
There is no problem with face maps (taking ``containment lemma'' for granted),
and the degeneracy maps are admissible in the sense 
that if $\bar{\Gamma }$ is the closure of the graph of the morphism in $\bar{\Delta }^q\times \bar{\Delta }^p$, then the inequality of Cartier divisors on $\bar{\Gamma }^N$:
\begin{equation*}
\bar{\Delta }^q_\infty |_{\bar{\Gamma }^N}\ge 
\bar{\Delta }^p_\infty |_{\bar{\Gamma }^N}
 \end{equation*}
holds.
%of \cite{KSY} which 
This is a sufficient condition for the preservation of the modulus condition,
and it can be shown that the admissibility may be checked after replacing the source by a surjective proper scheme over it equipped with the pull-back of the Cartier divisor, both by \cite[Lem.2.2]{KrishnaParkAdditive}.
We sketch how to verify the admissibility: By coordinate change it suffices to consider the projection $\bbA ^{n+1}=\Spec (k[x_1,\dots ,x_{n+1}])$ $\to $ $\bbA ^n=\Spec (k[x_1,\dots ,x_n])$ discarding the last factor. The ill-defined locus of the corresponding rational map between compactifications $\bbP ^{n+1}=\Proj (k[X_0,\dots ,X_{n+1}])$ $\dashrightarrow $ $\bbP ^n=\Proj (k[X_0,\dots ,X_{n}])$
is the closed subscheme $V(X_0,\dots ,X_n)$.
If we blow it up, the rational map becomes a morphism; the blow-up is covered by $(n+1)$-dimensional affine spaces
$\Spec (k[\frac{X_0}{X_j},\dots ,\frac{X_{n}}{X_j}, \frac{X_j}{X_{n+1}}])$ and other less important open subsets.
On this open subset, the pull-back of $\bar{\Delta }^{n+1}_\infty $
amounts to the divisor $\left( \frac{X_0}{X_j}\cdot \frac{X_j}{X_{n+1}} \right)$,
while that of $\bar{\Delta }^n_\infty $ amounts to $\left( \frac{X_0}{X_j} \right)$.
This verifies the admissibility.

%%%%%%%%%%%%%
%% Theorem %%
\begin{definition}
Let ${z}^i_{\mrm{simp}}(X|D,p)$ be the group of codimension $i$ cycles on $X\times \Delta ^p$ whose support satisfies the modulus condition and the face condition.
This gives us a simplicial abelian group $z^i_{\mrm{simp}}(X|D,\bullet )$.
\end{definition}
%% Theorem ends %%
%%%%%%%%%%%%%%%%%%

%%%%%%%%%%%%%%%%%%%%%%%%%%%%%%%
%%%%%%%%%%%%%%%%%%%%%%%%%%%%%%%
\subsection{Moving lemma}

The analogue of Theorem \ref{Th:full-statement} holds for the simplicial variant.
For this, the only new point is to establish an analogue of Proposition \ref{Prop:s(V)}
(see below for the analogue of the cheaper way as in Remark \ref{Rem:cheaper-way}).
We state and prove it now.

Let $\OF $ be a discrete valuation ring and $\pi \in \OF $ in the maximal ideal.
For a vector $\bm{v}\in \bbA ^n(\OF )$ and an integer $s\ge 0$, let $\Psi _{\bm{v},s}$ be the morphism
\begin{equation*}
\begin{array}{rcl}
\Psi _{\bm{v},s} \colon
\bbA ^n_{\OF }\times \bbA ^1&\to &\bbA ^n_{\OF } \\
(\bm{x},t)&\mapsto & (\bm{x}+t\pi ^s \bm{v}).
 \end{array}
\end{equation*}
We are interested in the process of pulling back a cycle in $\bbA ^n_\OF \times \Delta ^p$ by $\Psi _{\bm{v},s}\times \id _{\Delta ^p}$ and then by
one of the usual triangulation maps 
$\bbA ^n_\OF \times \Delta ^{p+1}\xrightarrow[\id _{\bbA ^n_\OF } \times \rho _i]{\cong } \bbA ^n_\OF \times (\bbA ^1\times \Delta ^p)$,
$i=1,\dots ,p+1$.

%%%%%%%%%%%%%
%% Theorem %%
\begin{proposition}
Given a closed subset $V$ of $\bbA ^n_\OF \times \Delta ^p$ satisfying the modulus condition with respect to the divisor $V(\pi )$,
There is an integer $s(V)\ge 0$ such that the pull-back 
\begin{equation*}
(\id _{\bbA ^n_\OF } \times \rho _i)\inv (\Psi _{\bm{v},s}\times \id _{\Delta ^p})\inv (V)
 \end{equation*}
satisfies the modulus condition in $\bbA ^n_\OF \times \Delta ^{p+1}$ for every $\bm{v}$ and $s\ge s(V)$.

\end{proposition}
%% Theorem ends %%
%%%%%%%%%%%%%%%%%%

%%% Proof begins %%%
\begin{proof}
By coordinate change, it suffices to consider the following morphism:
$\Delta ^{p+1}\cong \Spec (\bbZ [t_1,\dots ,t_{p+1}])$, $\Delta ^p\cong \Spec (\bbZ [t_1,\dots ,t_p])$ and the map is given by
\begin{equation*}
\begin{array}{ccc}
\bbA ^n_{\OF }\times \Delta ^{p+1}
&\to &\bbA ^n_{\OF }\times \Delta ^{p}
\\
(\bm{x},t_1,\dots ,t_{p+1})
&\mapsto & (\bm{x}+t_{p+1}\pi ^s \bm{v}, t_1,\dots ,t_p)
 \end{array}
 \end{equation*}
which we denote by $\Psi $ anew.
Denote the compactifications by
\begin{equation*}
\bar{\Delta }^{p+1}:=\Proj (\bbZ [U,T_1,\dots ,T_{p+1}]) \text{ and }
\bar{\Delta }^{p}:= \Proj (\bbZ [U,T_1,\dots ,T_{p}]).
 \end{equation*}

Modulus condition can be checked after blowing up the source, just like the admissibility of a morphism, by \cite[Lem.2.2]{KrishnaParkAdditive}.
We blow up $\bar{\Delta }^{p+1} $ by the ideal
$(U,T_1,\dots ,T_p)$ 
so that the rational map $\bar{\Delta }^{p+1}\dashrightarrow  \bar{\Delta }^p$ becomes well-defined. 

The scheme $\bar{\Delta }^{p+1}$ is covered by the open subsets where one of $U$, $T_1$, $ \dots $, $T_{p+1}$ is invertible.
Over the subset $\{ U\neq 0 \} =\Delta ^{p+1}$, there is nothing to check because one easily sees that the pull-back of $\bar{V}$ does not meet the pull-back of $D$ in this region.
The blow-up of the other open sets are covered by 
$(p+1)$-dimensional affine spaces 
$\tilde{\Delta }^{p+1}_{ij}$
having
$\frac{U}{T_j}$, $\frac{T_1}{T_j}$, $\dots $, 
$\widehat{\frac{T_i}{T_j}}$, 
$\dots $, 
$\frac{T_{p+1}}{T_j}$, $\frac{T_j}{T_i}$
as coordinates,
where $i=1, \dots , p+1$ and $j=1,\dots ,p $ with $i\neq j$.

Now we have a morphism induced by $\Psi $:
\begin{align*}
\Psi \colon &\bbA ^n_\OF \times 
\left( \tilde{\Delta }^{p+1}_{ij}
%\Spec (\bbZ \left[ \frac{U}{T_j}, \frac{T_1}{T_j},\dots ,\widehat{\frac{T_i}{T_j}}, \dots ,\frac{T_{p+1}}{T_j}, \frac{T_j}{T_i}
%\right] )
 \setminus \left\{  \frac{U}{T_j}=0 \right\} \right) \\
\to &\bbA ^n_\OF \times \Spec (\bbZ \left[ \frac{U}{T_i}, \frac{T_1}{T_i}, \dots ,\frac{T_p}{T_i}\right] )
\quad
\left( \overset{\text{open}}{\subset } \bbA ^n_\OF \times \bar{\Delta }^p\right) .
 \end{align*}
whose first component is given by $\bm{x}+\left( \frac{T_{p+1}}{T_j} / \frac{U}{T_j}\right) \pi ^s \bm{v} $,
and the second by the ring homomorphism
$\frac{(-)}{T_i}\mapsto \frac{(-)}{T_j}\cdot \frac{T_j}{T_i} $ where $(-)=U,\ T_1,\ \dots ,\ T_p$.
Let us denote by $\bar{V}_i$ the intersection of $\bar{V}$ and the target.
We now have to pull it back, take its closure in $\bbA ^n_\OF \times \tilde{\Delta }^{p+1}_{ij}$ and verify the inequality of Cartier divisors on its normalization required for the modulus condition.

Let $\{ f_\lambda (\bm{x}, \frac{U}{T_i}, \frac{T_1}{T_i}, \dots ,\frac{T_p}{T_i})\} _{\lambda }$ be finitely many functions defining $\bar{V}_i $ in the target.
Denote by $\deg _{\bm{x}}(f_\lambda )$ the total degree in the variables $\bm{x}$.
Then the defining ideal for $\left( \Psi \inv \bar{V}_i\right) ^-$ in $\bbA ^n_\OF \times \tilde{\Delta }^{p+1}_{ij}$ contains functions 
\begin{equation*}
\varphi _\lambda :=
\left( \frac{U}{T_j}\right) ^{\deg _{\bm{x}}(f_\lambda )}
\cdot
f_\lambda \left( \bm{x}+\left( \frac{T_{p+1}}{T_j} / \frac{U}{T_j}\right) \pi ^s \bm{v},\ 
\frac{U}{T_j}\frac{T_j}{T_i},\ \dots ,\ 
\frac{T_p}{T_j}\frac{T_j}{T_i}\right) .
 \end{equation*}
Expanding it with respect to the first entry,
it has the form
\begin{align}\label{phi-lambda}
\varphi_\lambda  = &\left( \frac{U}{T_j}\right) ^{\deg _{\bm{x}}(f_\lambda )}
\cdot
f_\lambda \left( \bm{x},\ 
\frac{U}{T_j}\frac{T_j}{T_i},\ 
\dots ,\ 
\frac{T_p}{T_j}\frac{T_j}{T_i}\right) 
\\
& + \pi ^s g_\lambda \left( \bm{x},\ \frac{U}{T_j},\
\dots ,\ \widehat{\frac{T_i}{T_j}},\ 
\dots ,\ \frac{T_{p+1}}{T_j},\ \frac{T_j}{T_i}\right) \notag
 \end{align}
where $g_\lambda $ is a polynomial $\in \OF [\bm{x}, \frac{U}{T_j},\dots ,\widehat{\frac{T_i}{T_j}},\dots ,\frac{T_{p+1}}{T_j},\frac{T_j}{T_i}]$.

On the other hand, since ${V}$ is assumed to satisfy the modulus condition, on $\bbA ^n_\OF \times \Spec (\bbZ \left[ \frac{U}{T_i}, \frac{T_1}{T_i}, \dots ,\frac{T_p}{T_i}\right] )$ we have a relation of the form 
\begin{equation*}
E_{i}( \frac{U}{T_i},\pi )
=\sum _{\lambda } b_\lambda f_\lambda 
 \end{equation*}
where $E_i(\alpha ,\beta )E_i[\bm{x},\frac{U}{T_i},\dots \frac{T_p}{T_i}](\alpha ,\beta )\in \left( \OF [\bm{x},\frac{U}{T_i},\dots ,\frac{T_p}{T_i}]\right) [\alpha ,\beta ]$ is a polynomial monic in $\alpha $ and homogeneous in $\alpha ,\beta $,
and $b_\lambda \in \OF [\bm{x},\frac{U}{T_i},\dots ,\frac{T_p}{T_i}] $.

We apply the $\OF [\bm{x}]$-algebra homomorphism $\frac{(-)}{T_i}\mapsto \frac{(-)}{T_j}\frac{T_j}{T_i}$ to get
\begin{align}\label{Eq:E-tilde}
\tilde{E}_i( \frac{U}{T_j}\frac{T_j}{T_i},\pi )
&:=
E_{i}[\bm{x},\frac{U}{T_j}\frac{T_j}{T_i},\dots \frac{T_p}{T_j}\frac{T_j}{T_i}]( \frac{U}{T_j}\frac{T_j}{T_i},\pi )
\\
&=\sum _{\lambda } b_\lambda (\bm{x},\
\frac{U}{T_j}\frac{T_j}{T_i},\ 
\dots ,\ 
\frac{T_p}{T_j}\frac{T_j}{T_i})
f_\lambda (\bm{x},\
\frac{U}{T_j}\frac{T_j}{T_i},\ 
\dots ,\ 
\frac{T_p}{T_j}\frac{T_j}{T_i})
\notag \end{align}
Now we assume $s\ge \deg _{\bm{x}}(f_\lambda )$ for all $\lambda $, and take
$\sum _\lambda b_\lambda \left( \frac{U}{T_j}\right) ^{s-\deg _{\bm{x}}(f_\lambda )} \cdot (-) $ of formulas \eqref{phi-lambda}; by formula \eqref{Eq:E-tilde}, we get 
\begin{align*}
(\text{Something } &\text{in the defining ideal for }\left( \Psi \inv \bar{V}_i\right) ^-)
\\
&=
\left( \frac{U}{T_j}\right) ^s\tilde{E}_i(\frac{U}{T_j}\frac{T_j}{T_i},\pi )
+  \pi ^s \left( \frac{U}{T_j}\right) ^{s
-\max _{\lambda }\{ \deg _{\bm{x}}(f_\lambda ) \} }
\tilde{g} \end{align*}
with $\tilde{g}\in \OF [\bm{x}, \frac{U}{T_j},\dots ,\widehat{\frac{T_i}{T_j}},\dots ,\frac{T_{p+1}}{T_j},\frac{T_j}{T_i}]$.
We multiply both sides by $\left( \frac{T_j}{T_i}\right) ^s$;
then if $2s-\max _{\lambda }\{ \deg _{\bm{x}}(f_\lambda )\} \ge s+\deg ({E}_i)$ where $\deg ({E}_i)$ denotes the homogeneous degree in $\alpha ,\beta $,
the right hand side becomes the form $\tilde{E}_i'(\frac{U}{T_j}\frac{T_j}{T_i},\pi ) $
with $\tilde{E}'_i (\alpha ,\beta )$ yet another homogeneous polynomial of degree $s+\deg (\tilde{E}_i)$.

Therefore if we let $s(V)$ the maximum of numbers $\deg (E_i) + \deg _{\bm{x}}(f_\lambda )$
(the index set for $\lambda $ varies with $i$),
our assertion holds.
\end{proof}
%%% Proof ends %%%

\begin{remark}\label{Rem:cheaper-simplicial}
As is the case with the cubical version,
the moving lemma for the simplicial version can be proved using $\Psi _{\bm{v},1}$ and push-forward by high enough powers $\bbA ^1\to \bbA ^1$.
For this, it is convenient to use a kind of sup-modulus condition on $\bbA ^n_\OF \times \bbA ^1\times \Delta ^p $ where a closed set $V $
is said to satisfy the modulus condition if locally on the normalization $\bar{V}^N$ of the closure of $V$ in $\bbA ^n_\OF \times \bbP ^1\times \bar{\Delta }^p$, either (where $D$ is the divisor defined by the function $\pi $)
\begin{equation*}
D|_{\bar{V}^N}\le F_1|_{\bar{V}^N}
\quad \text{ or }\quad
D|_{\bar{V}^N}\le \bar{\Delta }^p_\infty |_{\bar{V}^N}
 \end{equation*}
holds; then it can be checked that the modulus condition is preserved by pull-back along triangulation maps
$\Delta ^{p+1}\xrightarrow{\cong } \bbA ^1\times \Delta ^p$.

Let me mention the issue of the equivalence of cubical and simplicial version.
Miyazaki \cite{MiyazakiPrivate} has proved that the cubical and simplicial higher Chow groups with modulus are isomorphic if we go pro over the multiples of $D\subset X$. 
In the proof, he considers a double complex $z^i(X;p,q)$ of cycles in $X\times \Delta ^p\times \square ^q$, and the two associated spectral sequences,
just as in \cite{BlochCubical}.

One of the two spectral sequences degenerates thanks to the cube invariance of the cubical version (\cite{Miyazaki}, extended from the minus cube to pairs $(\bar{\Delta } ^p, -\bar{\Delta }^p_\infty )$). 
He observes that the simplicial version satisfies the cube invariance at least in pro;
this is because the triangulation maps
$\Delta ^{p+1}\xrightarrow{\cong }\bbA ^1\times \Delta ^p$ are not admissible but they are, up to doubling the divisor on the source.
This implies that the other spectral sequence degenerates at least in pro, which enables him to establish a parallel of the argument in \cite{BlochCubical}.

\end{remark}

\bibliographystyle{alpha}%achicago = Chicago style
\bibliography{WataruBib}

\end{document}